\documentclass{article}
\usepackage[utf8]{inputenc} 
\usepackage{enumerate}
\usepackage[dvipsnames]{xcolor}

\usepackage{amsmath,indentfirst,stmaryrd,qsymbols,graphicx,amsfonts,amsthm,ulem,dsfont,cleveref,footmisc,scalerel,stackengine,xspace,amsmath}
\usepackage{fourier}
\usepackage{subcaption}
\usepackage[margin=2cm]{geometry}
\def \d{{\mathrm d}}

\def \Vol{{\sf Vol}}

\def \Floor{{\sf Floor}}
\def \SubPrism{{\sf SubPrism}}

\def \CH{{\sf CH}}
\def \CCSF{{\sf CCSF}}
\def \Mountain{{\sf Mountain}}

\def \1{\mathds{1}}

\def \AP{{\sf AP}}

\definecolor{amber}{rgb}{1.0, 0.75, 0.0}

\def \Layer{{\sf Layer}}

\def \N{\mathbb{N}}

\def \P{\mathbb{P}}
\def \R{\mathbb{R}}

\def \Z{\mathbb{Z}}

\def \app#1#2#3#4#5{\begin{array}{rccl} #1:&#2&\longrightarrow&#3\\ &#4&\longmapsto&#5\end{array}}

\def \bM{\begin{bmatrix}}

	\def \ba{\begin{align}}
		\def \ba{{\bf a}}

		\def \ben{\begin{eqnarray}}
			\def \beq{\begin{equation}}
				\def \be{\begin{eqnarray*}}

					\def \bia{\begin{itemize}\compact \setcounter{d}{0}}
						
						\def \bir{\begin{itemize}\compact \setcounter{c}{0}}
							\def \bis{\begin{itemize}\compact }
								\def \bi{\begin{itemize}\compact \setcounter{b}{0}}
									
									\def \bls{{\tiny $\blacksquare$ }}
									
									\def \bpar#1{\left\{\begin{array}{#1} }

										\def \build#1#2#3{\mathrel{\mathop{\kern 0pt#1}\limits_{#2}^{#3}}}

										\def \captionn#1{\begin{center}\begin{minipage}{15cm}\sf\caption{\small #1}\e
													d{minipage}\end{center}}

											\def \dd#1{\frac{\partial}{\partial #1}}
											\def \dd{\xrightarrow[n]{(d)}}

											\def \eM{\end{bmatrix}}
										\def \ea{\end{align}}
									\def \een{\end{eqnarray}}
								\def \ee{\end{eqnarray*}}
							\def \eia{\end{itemize}\vspace{-2em}~}
						\def \eir{\end{itemize}\vspace{-2em}~}
					\def \eis{\end{itemize}\vspace{-2em}~}
				\def \ei{\end{itemize}\vspace{-2em}~}
			
			\def \epar { \end{array}\right.}

		\def \eq{\end{equation}}
	\def \eref#1{(\ref{#1})}

	\def \imp{\Rightarrow}
	
	\def \Dom#1{{\sf Dom}(#1)}
	
	\def \l{\left}

	\def \r{\right}

	\def \sous#1#2{\mathrel{\mathop{\kern 0pt#1}\limits_{#2}}}
	\def \sur#1#2{\mathrel{\mathop{\kern 0pt#1}\limits^{#2}}}

	\def \Conca{{\sf Conca}}
	
	\def \abso#1{\vert#1\vert}

	\def\bma{ \begin{bmatrix}}
		\def\ema{\end{bmatrix}}
	\def\cro#1{\llbracket#1\rrbracket}

	\newcommand{\E}{\mathbb{E}}
	
	\newcommand{\ie}{\textit{i.e. }}

	\newcommand{\compact}{ \topsep0pt   \itemsep=0pt   \partopsep=0pt   \parsep=0pt}

	\newcounter{b}
	\newcounter{c}
	\newcounter{d}
	\newqsymbol{`B}{\mathcal{B}}
	\newqsymbol{`E}{\mathbb{E}}
	\newqsymbol{`I}{\mathbb{I}}
	\newqsymbol{`N}{\mathbb{N}}
	\newqsymbol{`O}{\Omega}
	\newqsymbol{`P}{\mathbb{P}}
	\newqsymbol{`Q}{\mathbb{Q}}
	\newqsymbol{`R}{\mathbb{R}}
	\newqsymbol{`W}{\mathbb{W}}
	\newqsymbol{`Z}{\mathbb{Z}}
	\newqsymbol{`a}{\alpha}
	\newqsymbol{`e}{\varepsilon}
	\newqsymbol{`o}{\omega}
	\newqsymbol{`t}{\tau}
	\newqsymbol{`w}{{\cal W}}
	\newtheorem{lem}{Lemma}[section]

	\newtheorem{conj}[lem]{Conjecture}
	
	\newtheorem{defi}[lem]{Definition}

	\newtheorem{pro}[lem]{Proposition}
	\newtheorem{rem}[lem]{Remark}
	\newtheorem{theo}[lem]{Theorem}

\begin{document}
		\author{Jean-François Marckert\footnote{Univ. Bordeaux, CNRS, Bordeaux INP, LaBRI, UMR 5800, F-33400 Talence, France}\and Ludovic Morin\footnotemark[2]}
	\title{The Sylvester question in $\R^d$: convex sets with a flat floor\thanks{The authors were partially supported by the  ANR projects   CartesEtPlus (ANR-23-CE48-0018) and 3DMaps (ANR-20-CE48-0018)}}
	\date{}
	\maketitle
	
	\abstract{ Pick $n$ independent random points $U_1,\ldots,U_n$ uniformly in a compact convex set $K$ of $\R^d$ with volume 1, and let $P^{(d)}_K(n)$ be the probability that these points are in convex position. The Sylvester conjecture in $\R^d$ is that $\min_K P^{(d)}_K(d+2)$ is achieved by $d$-dimensional simplices $K$ (only). \par 
		In this paper, we focus on a companion model, already studied in the $2D$ case, which we define in any dimension $d$: we say that $K$ has $F$ as a flat floor, if $F$ is a subset of $K$, contained in a hyperplane $P$, such that $K$ lies in one of the half-spaces defined by $P$. \par
		We define $Q_K^F(n)$ as the probability that $U_1,\cdots,U_n$ together with $F$ are in convex position (i.e., the $U_i$ are on the boundary of the convex hull ${\sf CH}(\{U_1,\cdots,U_n\}\cup F\})$). We prove that, for all fixed $F$, 
		$K\mapsto Q_K^F(2)$ reaches its minimum on the "mountains" with floor $F$ (a mountain is the convex hull of $F$ together with an additional vertex), while the maximum is not reached, but $K\mapsto Q_K^F(2)$ has values arbitrary close to 1. If the optimisation is done on the set of $K$ contained in $F\times[0,d]$ (the "subprism case"), then the minimum is also reached by the mountains, and the maximum by the "prism" $F\times[0,1]$. The probability $P^{(d)}_K(d+2)$ can be expressed in terms of the expected volume of ${\sf CH}(\{U_1,\dots,U_{d+1}\})$ and analogously, $Q_K^F{(2)}$ relies on the expected volume of ${\sf CH}(\{U_1\}\cup F\})$, so that our optimization result can be seen as a proof of the Sylvester's problem in the floor case. \par
		In $2D$, where $F$ can essentially be seen as the segment $[0,1],$ we give a general decomposition formula for $Q_K^F(n)$, which can be used to compute several formulas and bounds for different $K$. In 3D, we give some bounds for $Q_K^F(n)$ for various floors $F$ and special cases of $K$.}

	\section{Introduction and main results}
	Let $d\geq2$ be a fixed integer. Throughout the paper, the space $\R^d$ is equipped with an orthonormal basis $e_1,\cdots,e_d$, and its associated coordinate system. 
	
	For $n\geq 1$, we say that an $n$-tuple of points $u[n]=(u_1,\ldots,u_n)\in(\R^d)^n$ is in convex position if $\{u_1,\cdots,u_n\}$ forms the vertex set of a convex polytope. For any compact convex set $K$ in $\R^d$ with non-empty interior and for any $n\in\N$, we let  $\mathbb{U}_K$ be the uniform distribution on $K$, and let $\mathbb{U}_K(n)$ be the product measure of $n$ iid random variables $U[n]=(U_1,\cdots,U_n)$ with common distribution $\mathbb{U}_K$.

	Finally, we let $P^{(d)}_K(n)$ be the probability that an $n$-tuple of points $U[n]$ with distribution $\mathbb{U}_K(n)$ is in convex position.
	
	The very first problem in the study of $P^{(d)}_K(n)$ was posed by Sylvester \cite{sylvester} in 1864 (see Pfiefer \cite{pfiefer} for a history of this problem) and consisted in determining the probability that 4 random points in the plane, are in convex position. This problem evolved over time into computing the convex sets $K$ that minimize or maximize the probability $P^{(2)}_K(4)$. Blaschke \cite{blaschke1917affine} provided the answer to this question in 1917, by proving that the lower bound is reached when $K=\triangle^2$ is a triangle, and the upper bound when $K=\bigcirc^2$ is a disk or an ellipse, namely
	\[\frac{2}{3}=P^{(2)}_{\triangle^2}(4)\leq P^{(2)}_K(4)\leq P^{(2)}_{\bigcirc^2}(4)=1-\frac{35}{12\pi^2}.\]
	In larger dimensions, an analogous problem exists. Since $d+1$ random points are in convex position with probability 1, 
	the question is whether we may optimize the quantity $P^{(d)}_K(d+2)$ ?
	As a generalization of Blaschke's proof, Groemer \cite{groemer1974mean} proved that the upper bound $P^{(d)}_K(d+2)\leq P^{(d)}_{\bigcirc^d}(d+2)$ holds true whereas the following $d$-Sylvester's question 
	\[P^{(d)}_{\triangle^d}(d+2)\leq P^{(d)}_K(d+2)\leq P^{(d)}_{\bigcirc^d}(d+2),\]
	($\triangle^d$ and $\bigcirc^d$ denote respectively a simplex and an ellipsoid of volume 1) remains a conjecture due to the lower bound. To prove the lower bound, Blaschke uses the fact that if $K$ is a convex set, and $S(K)$ is obtained from $K$ by a ``shaking'' with respect to a given line, then $P^{(2)}_{S(K)}\leq P^{(2)}_K$ (this eventually allows to conclude, because, for each $K$, there exists a sequence of shakings $(S_i)$ such that $S_i(...S_1(K))$ converges to a triangle). A shaking with respect to an hyperplane is well defined in higher dimension, but the property  $P^{(d)}_{S(K)}\leq P^{(d)}_K$ does not hold in all generality.
	
	The quantity $P^{(d)}_K(d+2)$ satisfies
	\ben\label{eq:ef} P^{(d)}_K(d+2)= 1- (d+2)\E\l( \Vol_d({\sf CH}(\{U_1,\cdots,U_{d+1}\}))\r),\een 
	(where ${\sf CH}(\{U_1,\cdots,U_{d+1}\})$ is the convex hull of $\{U_1,\cdots,U_{d+1}\}$)
	and this formula may be understood as the fact that $U_1,\cdots,U_{d+2}$ are not in convex position if and only if one of the $U_i$ is in the convex hull of the others.
	Hence, the problem of optimizing $K\mapsto P^{(d)}_K(d+2)$ is equivalent to the problem of optimizing $K\mapsto -\E( \Vol_d({\sf CH}(\{U_1,\cdots,U_{d+1}\})))$.
	
	This problem seems to be out of reach for the moment, and then, we will turn our attention to another problem, which has a common flavor for many reasons, which should become clear gradually.
	
\begin{figure}[h!] 
		\centerline{\includegraphics{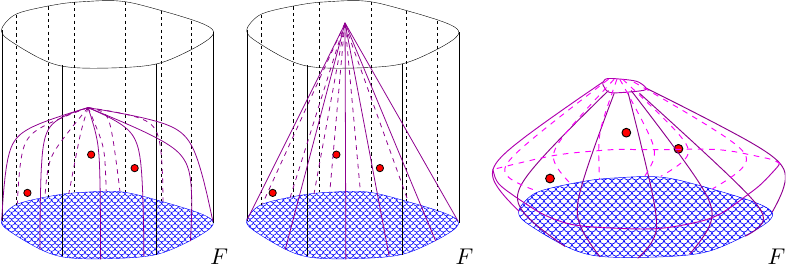}}
		\begin{subfigure}[t]{0.3\textwidth}
			\centering
			\caption{A compact convex set with floor $F$. It is a subprism.}
		\end{subfigure}%
		~ 
		\begin{subfigure}[t]{0.28\textwidth}
			\centering
			\caption{A mountain having same floor $F$ as the left one. It is also a subprism.}
		\end{subfigure}
		~ 
		\begin{subfigure}[t]{0.35\textwidth}
			\centering
			\caption{This third compact convex set is not a subprism but is an element of $\CCSF(F)$.}
		\end{subfigure}
		\caption{\label{fig:twofam}{Examples of compact convex sets of $\R^3$ with same floor $F$.}}
	\end{figure}

	\begin{defi} Let $K$ be a compact convex set of $\R^d$, for some $d \geq 2$, and $F$ be a subset of $\R^{d-1}$.  
		We say that $F$ is the floor of $K$ if  
		$K$ is contained in the half space  $\{x_d\geq 0\}$, and if its intersection with the hyperplane $x_d=0$ is $F\times \{0\}$.
	\end{defi}
	Hence, the floor of $K$ exists if $K$ intersects the hyperplane $\R^{d-1}\times\{0\}$ and lies above it. By definition, a floor $F$ is always a compact convex set of $\R^{d-1}$.
	For the sequel, we will need to distinguish \textbf{two families of convex sets} (see Fig. \ref{fig:twofam}):
	\begin{defi}
		\bls We say that a convex set $K$ with floor $F$ is a \textbf{sub-prism} if $K$ is a subset of $F \times[0,H]$ for some $H\geq 0$. \par
		\bls We denote by $\SubPrism(F)$ the set of convex sets $K$, that are sub-prisms with floor  $F$, and such that $\Vol_d(K)=1$.  
	\end{defi}
	In the following, we will say that  $\Floor(K) \times[0,H]$ is a prism.
	
	In geometry, the word "prism" is often used in 3D to identify a polyhedron that is equal to the convex hull of two (non-co-planar) polygons, image of each other by a translation. In this case, each face other than $F$ is still a polygon. It will not be the case here, since we do not necessarily assume that $F$ is a polytope (it can be smooth, for example).

	The second family we introduce is: 
	\begin{defi} We denote by $\CCSF(F)$ the set of Compact Convex Sets $K$ with  Floor $F$, and such that $\Vol_d(K)=1$. 
	\end{defi}
	
	\paragraph{Convention:} Throughout the paper, $F$ will be a compact convex set of $\R^{d-1}$ such that 
	\ben{\sf Vol}_{d-1}(F)=1.\een
	Since affine maps preserve both convexity and the uniform distribution, this is not a real restriction.

	\begin{figure}[h!]
		\centerline{\includegraphics[scale=0.9]{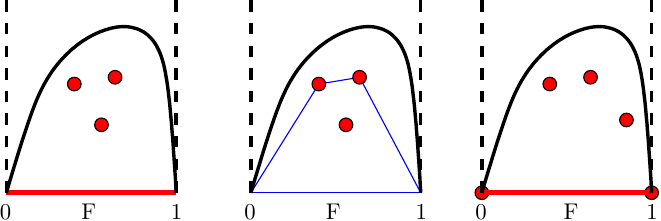}}
		\caption{In 2D, a subprism convex set with floor $[0,1]\times\{0\}$; in the first picture the 3 points are in convex position, but they are not in convex position together with $F$, since, as shown in the second picture, the convex hull of the union of $F$ and of the three points does not contain all these objects on its border. In the third picture, they are in convex position together with $F$.}
	\end{figure}	
	
	\paragraph{A modified Sylvester question in the fixed floor case.}
	
	Let us fix a floor $F$, which will always be assumed to be a compact convex set  (of volume $1$) in $\R^{d-1}$. For any $K\in \CCSF(F)$, denote by $Q_K(2)$ the probability that the convex hull of the set formed by the floor of $K$ (that is, the set $F$), and two random points taken independently and uniformly in $K$ contains these two points on its boundary.
	
	We will motivate this question further  in Section \ref{sec:motivation}, but for the moment, observe that the problem of optimizing $K\mapsto Q_K(2)$ has the same flavor as that of $K\mapsto P^{(2)}_K(4)$ due to an analogue of the Efron formula \cite{efronformula}; indeed, by the same reasoning as for \eref{eq:ef} we have:
	\ben\label{eq:gefron} Q_K(2) = 1- 2\, \E\Big({\Vol_d}(\CH(\{U\} \cup F)\Big),\een
	where $U$ is uniform in $K$.
	The two quantities $P^{(d)}_K(d+2)$ and $Q_K(2)$	are the first "non trivial quantities" in both cases.

	We then introduce a geometrical object:
	\begin{defi} 
		For a point $z$ with positive last coordinate $z_d\geq 0$, we call mountain with apex $z$, the compact convex set  
		\[\Mountain_F(z)=\CH(\{z\} \cup F).\]
		
	\end{defi}
	Of course, $\Mountain_F(z)$ has floor $F$, and since $\Vol_{d-1}(F)=1$,
	\ben{\sf Vol}_d(\Mountain_F(z))= z_d  / d .\een 
	Two of the main contributions of this paper are the following theorems: 
	
	\begin{theo}\label{theo:Prism}[Optimization in the set $\SubPrism(F)$]. Let $d\geq 2$, $F$ be a compact convex set in $\R^{d-1}$ such that $\Vol_{d-1}(F)=1$. 
		For all sub-prism $K\in \SubPrism(F)$, 
		\[ Q_{{\sf UnitMountain}}(2) \leq Q_K(2) \leq Q_{\sf UnitPrism}(2),\]
		where the ${\sf UnitPrism}$ is the prism $F\times[0, 1]$  and a ${\sf UnitMountain}$ identifies any mountain with floor $F$ and volume 1.\par
		Moreover, 
		\ben\label{eq:res}
		Q_{{\sf UnitMountain}}(2) &=& 
		1-{2}/{(d+1)}\\
		Q_{\sf UnitPrism}(2)&=& 
		1-1/d,\een both of which do not depend on $F$. 
	\end{theo}
	Again, since invertible affine maps preserve convexity and the uniform measure, this result applies to "non vertical prisms".
	\begin{rem}
		The bound $Q_K(2) \leq Q_{\sf UnitPrism}(2)$ holds even if $K$ is not convex. The proof we will give entails that if $K$ can be written as $\cup_{z\in F} \{z\}\times [0,G(z)]$ where $G(z)$ is the height at $z$, then under the only hypothesis that $G$ is measurable and satisfies $\int_F G(z)dz=1$, then in this case  $Q_K(2) \leq Q_{\sf UnitPrism}(2)$.
	\end{rem}
	\begin{theo}\label{theo:CCSF}[Optimization in the set $\CCSF(F)$]. Let $d\geq 2$, $F$ be a compact convex set in $\R^{d-1}$ such that $\Vol_{d-1}(F)=1$. 
		For all  $K\in \CCSF(F)$, 
		\[ Q_{{\sf UnitMountain}}(2) \leq Q_K(2) < 1 ,\]
		and there exists an element of $\CCSF(F)$ such that $Q_K(2)\geq 1- \alpha$, for all $\alpha>0$, and none such that $Q_K(2)=1$.
	\end{theo}
	
	\normalsize
	
	\subsection{Related results in 2D} Under the quantity $P^{(d)}_K(n)$ lie many other questions. Can we compute $P^{(d)}_K(n)$ for some specific sets $K$? Can we optimize $P^{(d)}_K(n)$ on $K$ for other values of $n$? Is there a limit shape theorem for a $n$-tuple of points $U[n]$ that are $\mathbb{U}_K(n)$-distributed and conditioned to be in convex position ?
	
	To date, the overwhelming majority of known answers to these problems have been given when $d=2$, but even there the questions are still largely open.

	For example, a natural generalization of Sylvester's problem  when $d=2$ would be 
	\[P^{(2)}_\triangle(n)\leq P^{(2)}_K(n)\leq P^{(2)}_\bigcirc(n),\] for any compact convex set $K$ of $\R^2$. The case $n=5$ was proved by Marckert and Rahmani \cite{marckert:hal-02913348}, but $n\geq 6$ is still a conjecture (but one finds in \cite{marckert:hal-02913348} an explicit formula to compute $P^{(2)}_K(n)$). 
	
	As for the computation of $P^{(2)}_K(n)$, several results are known. If $K$ is a parallelogram or a triangle, Valtr \cite{Valtr1995, valtr1996probability} computed two exact formulas for all $n\geq3$:
	\[{P}^{(2)}_\Box(n)=\frac{1}{(n!)^2}{2n-2\choose n-1}^2\quad\text{ and }\quad{P}^{(2)}_\triangle(n)=\frac{2^n (3n-3)!}{(2n)!((n-1)!)^3}.\]
	
	\noindent\textbf{In dimension $d=2$, a fixed floor $F$ }can only be a nonzero segment. Thus, $Q_K(n)$ is the probability that $n$ points in $K$ are in convex position, together with both extremities of $F$.  In fact, 
	\begin{figure}[htbp]
		\centerline{	\includegraphics{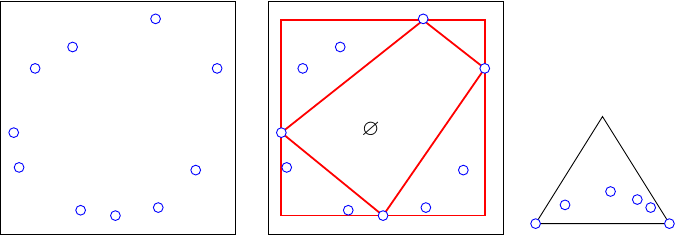}}
		\caption{\label{fig:2d-dec} Buchta's method \cite{buchta2012} to compute the probability that the $n$ points $U_1,\cdots,U_n$ are in convex position in a parallelogram. Push the sides of the initial parallelogram inside, and stop them when they hit a point. In the end, each of the four sides contains a point (the 4 points are not necessarily different). This forms 2, 3 or 4 triangles, all of which having 2 marked points at their boundaries. The $U_i$ are in convex position, if and only if each of the other points belong to a triangle, and if  all the points inside a triangle together with the two marked points are in convex position.}
	\end{figure}
	an important tool for these results relates to the $2D$-floor case. Indeed, a classical method to compute $P_K(n)$ consists "in shrinking $K$", e.g., by pushing the boundaries inwards (as in Buchta \cite{buchta_2006}, Marckert \cite{marckert2017probability}, Morin \cite{morin2024,morin2024bis}, see also Fig. \ref{fig:2d-dec}), until some points among the $U_i$ are "hit". This defines several regions, which either must not contain any points, or must bound a compact convex set containing a subset of the $U_i$ that must be in convex position together with the two extremities of a segment (these extremities being some $U_i$ discovered when pushing the boundaries). Up to an affine map, this latter set is a compact convex set with a floor.
	In \cite{buchta_2006} and \cite{morin2024,morin2024bis}, the zones are triangles, while in \cite{marckert2017probability} they are the interior of some arcs of disk parameterized by an angle $\theta$.  
	
	Some important asymptotic results also involve this type of tool, such as Bárány's formula \cite{barany2}, who proved that for any compact convex set $K\subset \R^2$, 
	\[\lim_{n\to+\infty} n^2\left({P}^{(2)}_{K}(n)\right)^{\frac{1}{n}}=\frac{1}{4}e^2\AP^*(K)^3,\]
	where $\AP^*(K)$ is the supremum of the affine perimeters of all convex sets $S\subset K.$ In the case where $K$ is a polygon, this logarithmic equivalent was refined into a true equivalent by  Morin \cite{morin2024,morin2024bis}.
	Both authors used the well-known case when $K$ is a triangle, due to Bárány, Röte, Steiger and Zhang \cite{baranybipointed}: 
	\begin{align}\label{eq:bipointee}
		t_n:=Q_\triangle(n)=\frac{2^n}{n!(n+1)!}.
	\end{align}
	Again, the notation "$Q$" refers to the floor case when $P$ refers to the classical one.
	
	Other questions then arise. Do we know other values of $Q_K(n)$? Can we optimize $Q_K(n)$ for values of $n\geq 3$ ? In this paper, we give a recursive formula for $Q_K(n)$ in the $2D$-case which allows to recover the triangle case (Theorem \ref{theo:tri}), to compute exact values in the case where $K$ is a square (denoted by $q_n$, see Proposition \ref{pro:q}) or a parabola (Theorem \ref{theo:para}) above the segment $F$ (denoted by $p_n$), namely
	\[p_n = \frac{12^{n+1}}{6 (2 n+2)!}\quad\text{ and }\quad q_n =\frac1{n!(n+1)!}{2n \choose n}.\]
	Of course, a rewording of Theorem \ref{theo:Prism} in the $2D$-case gives immediately that
	\[t_2=\frac{1}3\leq Q_K(2)\leq q_2=\frac12.\]
	However, an additional independent elementary proof of this result will be given in Section \ref{sec:secondproof}. 
	More generally, we believe that,
	\begin{conj}
		For all $n\geq 3$, and all $K\in \CCSF([0,1]),$
		\[t_n\leq Q_K(n)\leq q_n.\]
	\end{conj}
	
	Our belief in this conjecture relies on several facts, but the main one is the following: we made a random generator of convex sets $K$ with floors, on which we compared empirically $Q_K(n)$ to $t_n$ and $q_n$, for small values of $n$ (we tested hundreds of shapes).

	\subsection{Bounds and universal bounds in dimension $3$}
	
	When $d>2$, exact results for the computation of $Q_K(n)$ or $P^{(d)}_K(n)$ are rare. Most of the known results were obtained when $K=\bigcirc^d$, the Euclidean unit ball in dimension $d$. The value of $P^{(3)}_{\bigcirc^3}(5)$ (that is, 5 points in the $3D$-ball) was computed by Hostinsky \cite{Hostinsky}, before a general formula for $P_{\bigcirc^d}^{(d)}(d+2)$ was found by Kingman \cite{kingman}. Miles \cite{miles} proved that 
	$\lim_{d\to\infty}P_{\bigcirc^d}^{(d)}(d+3)=1,$ and later conjectured that for all $m>4$, $\lim_{d\to\infty}P_{\bigcirc^d}^{(d)}(d+m)=1,$ a conjecture proved by Buchta \cite{Buchta1986}.
	
	Even in dimension $3$, very little is known. The most classical and simplest results in dimension 2  have no analog in dimension 3 (however, 
	Buchta \& Reitzner \cite{BR2001} computed the expected volume of the convex hull of $n$ iid uniform points taken in the tetrahedron).
	
	The last section of this paper tries to remedy this problem by providing bounds in the floor case.
	In Theorem \ref{theo:dfreghrtu} , we will prove that for any convex set $F$ of volume $1$ and any unit mountain $M$ having floor $F$, the probability $Q_M(n)$ satisfies
	\[Q_M(n)\geq \frac{2^n}{n!}\prod_{j=1}^n \frac1{3j-1}.\]
	When $F$ is a triangle, such a mountain is a tetrahedron $T$.
	In this case,  we also obtained an upper bound (and an improved lower bound) in Theorems \ref{theo:upp} and \ref{theo:dqgrh}, and so we give two explicit sequences $(u_n)$ and $(\ell_n)$ such that 
	\ben \ell_n \leq Q_{T}^{(3)}(n) \leq u_n, \textrm{~~for~~} n\geq 0.\een

	\subsection{Motivation}\label{sec:motivation}
	There are several motivations for studying these "floor case questions".
	
	\bls As already mentioned, in $2D$, ``the floor case'' is already present in the literature (but the terminology ``floor case'' is not), since convex chains in a triangle (with two vertices of the triangle being elements of this chain) appear under some forms in many papers devoted to points in convex position: Bárány \& al. \cite{baranybipointed}, Valtr \cite{Valtr1995, valtr1996probability}, Buchta \cite{buchta_2006}, Morin \cite{morin2024bis,morin2024}.

	The computation of $P_K^{(d)}(d+2)$ becomes much more difficult as soon as $d\geq3$, and this is most likely due to the fact that there is no convenient decomposition suitable for a combinatorial approach. For example, if one draws $n$ iid uniform points $U_1,\cdots,U_n$ in the 3D cube, and if one pushes the faces of the cube inward until they hit a point in $U[n]$, then one gets a parallelepiped with marked points, and this structure does not seem to be helpful for further decomposition of the problem! (For example, if each of the 6 faces contains a distinct point, then, the complement of the convex hull of these points is still connected, and it seems difficult to make any progress from there).
	
	However, because of the success and simplicity of the study of the "floor case in 2D", it is natural to try to understand the floor case in 3D, and in the following dimensions.
	
	\bls Take a compact convex set $K$ of dimension $d$, and let $d+2$ iid uniform points $U_1,\cdots,U_{d+2}$ into $K$.
	For $1\leq m\leq d,$ let $S_{m}$ be the convex hull of $U_1,\cdots,U_m$, and $H_m$ be the affine space of dimension $m-1$  containing $S_m$. We have
	\be
	P_K^{(d)}(d+2)&=&1-(d+2) \E\Big( \Vol_{d}( S_{d+1})\Big)\\
	&=&1- \frac{d+2}d \E\Big( \Vol_{d-1}( S_{d}) \times  D(U_{d+1},H_{d}) \Big)
	\ee
	where $ D(U_{d+1},H_{d})$ is the Euclidean distance of  $U_{d+1}$ from the hyperplane $H_d$ (for more information, see e.g. Kingman \cite{kingman}, who computed $ P_K^{(d)}(d+2)$, where $K$ is the unit Euclidean ball in $\R^d$).
	
	To solve the  Sylvester problem, it is natural to try to understand the typical (joint) behavior of the two random variables $\Vol_{d-1}(S_{d})$ and $D(U_{d+1},H_{d})$, since the probability that $U_{d+2}$ belongs to $\CH(U[d+1])$ is $\Vol_{d-1}(S_{d}) \times D(U_{d+1},H_{d}) /d$.
	A small picture shows that $H_d$ cuts $K$ into two parts $P_1$ and $P_2$, and $U_{d+1}$ will fall into one of these parts, proportional to its $d$-volume. Hence, $H_d$ plays the role of a floor, and conditional on $H_d$ (and of the part $P_j$ in which $U_{d+1}$ is)  $D(U_{d+1},H_{d})$ has the same distribution "as the height of a random point in $P_j$". Trying to understand what kind of part maximizes the height expectation -- and we proved that they are mountains -- can be seen as a step  towards solving the Sylvester problem.

	\section{Optimization of \texorpdfstring{$Q_K(2)$}{Q_K(2)} in the multidimensional case}
	
	Consider an element $K$ taken in $\CCSF(F)$ or in $\SubPrism(F)$, and take $U$ a uniform random point inside $K$.
	Set $H=U_d$, the last coordinate of $U$, seen as a height, and let $Z=(U_1,\cdots,U_{d-1})$ the $d-1$ first coordinates. If $K$ is in $\SubPrism(F)$, then $Z\in F$ almost certainly, while this may not be the case for $K\in\CCSF(F)$. 
	
	Consider the "layer of $K$ at level $t$", defined by
	\[\Layer_K(t)= K \cap \l(\R^{d-1} \times \{t\}\r),\]
	and consider its volume
	\[ L_K(t):= \Vol_{d-1}(\Layer_K(t)).\]
	It is clear that $L_K$ is null outside the compact set $[0,m_K]$ where $m_K=\max\{ z_d~:~ z \in K\}$ is "its maximal height".   
	By standard calculus
	\[\Vol_d(K)=\int_0^{m_K} L_K(t) dt=1,\]
	and in fact, $L_K$ is the density of the random variable $H$. 
	\begin{lem} Let $K\in \CCSF(F)$ (which includes $\SubPrism(F)$),
		then 
		\ben\label{eq:QL} Q_K(2) = 1-2 \int_{0}^{+\infty} L_K(t) \frac{t}{d} dt.\een
	\end{lem}
	\begin{proof} By \eref{eq:gefron}, 
		\ben\label{eq:reph} Q_K(2)=1-2 \E\l(\Vol_{d}( \Mountain_F(Z,H))\r)=1-2\, \E({H}/d)= 1-2 \int_{0}^{+\infty} L_K(t) \frac{t}{d} dt.\een
	\end{proof}
	Hence, in the case of unitary mountains, 
	\ben L_{{\sf UnitMountain}}(t) =  \l(1-  t / d\r)^{d-1} \1_{[0,d]}(t),\een
	and then  $Q_{{\sf UnitMountain}}(2)= 1-2  \int_{0}^d \l( 1-t/d\r)^{d-1} \frac{t}{d} dt = 1-\frac{2}{d+1}$,
	as announced. In the case of ${\sf UnitPrism}$, $L_{{\sf UnitPrism}}(t)=\1_{[0,1]}(t)$, so that 
	$Q_{{\sf UnitPrism}}(2)= 1-2  \int_{0}^1 \frac{t}{d}  dt = 1-1/d$ as announced.

	We will first establish Theorem \ref{theo:Prism}, and then Theorem \ref{theo:CCSF}.

	\subsection{SubPrism case: proof of Theorem \ref{theo:Prism}}

	\subsubsection{Proof of Theorem \ref{theo:Prism}: upper bound}
	\label{sec:egther} 
	Each element $K\in \SubPrism(F)$ is characterized by its top function $G_K$ defined by
	\ben \app{G_K}{F}{[0,+\infty)}{z=(z_1,\cdots,z_{d-1})}{G_K(z):=\max\{y \in \R, (z_1,\cdots,z_{d-1},y) \in K\}}\een
	since 
	\ben\label{eq:Cara} K = \bigcup_{z \in F} \Big[ \big(z_0,\cdots,z_{d-1},0\big), \big(z_0,\cdots,z_{d-1},G_K(z)\big)\Big].\een 
	Thus, if $K\in \SubPrism(F)$, both $G_K$ and $K$ characterize each other. 
	Moreover, a function $G:F\to\R$ is the top function of a compact convex set $K$ in $\SubPrism(F)$ iff it is non-negative, concave, and satisfies $\int_F G(z)dz=1$.

	Recall that for $U=(Z,H)$ uniform in $K$, and by \eref{eq:reph}, $Q_K(2)=1-2\,\E(H/d)$. 
	The density of $Z$ is $G_K$ on $F$ while, given $Z$, $H$ is uniform on $[0,G_K(z)]$, so that $\E(H~|~Z=z)= G_K(z)/2$. Hence
	\ben \label{eq:relrhte}
	Q_K(2)&=&1-2\int_{F} G_K(z) \frac{G_K(z)}{2d} dz 
	=1-\int_{ F} \frac{G_K^2(z)}{d}    dz. \label{eq:cvx}
	\een 
	From here the proof of the upper bound is quasi-immediate:    
	\begin{align*}
		Q_{{\sf UnitPrism}}(2)-Q_K(2)=\int_{F} \frac{G_K^2(z)-1}{d}    dz 
		=\int_{ F} \frac{(G_K(z)-1)^2}{d}    dz + \underbrace{\int_F\frac{(2G_K(z)-2)}ddz}_{=0} \geq 0.
	\end{align*}

	\subsubsection{ Proofs of Theorems \ref{theo:Prism}  and of  \ref{theo:CCSF}: lower bound}

	The proof we give is valid in $\CCSF(F)$, and since the lower bound is reached for unitary mountains (for any two unitary mountains $M_1$ and $M_2$, $Q_{M_1}(2)=Q_{M_2}(2)$), some of which are in $\SubPrism(F)$, the proof that follows is valid in $\SubPrism(F)$ too.
	
	In the following, $K$ is an element of $\CCSF(F)$, and $M$ is a unitary mountain in $\CCSF(F)$. Using \eref{eq:reph}, it suffices to prove that
	\ben\label{eq:qdgr} \E(H_K) \leq \E(H_M) \een
	where $H_K$ and $H_M$ denote the respective random variable heights of a point $U_K$ and $U_M$ taken uniformly in $K$ and $M$. 
	For any non-negative random variable $X$, we have $\E(X)= \int_0^{+\infty} \P(X\geq t)\, dt$. 
	Since 
	\[\P(H_K\geq t)= \int_{t}^{\infty} L_K(s)\,ds=1-\int_{0}^{t} L_K(s)\,ds,\] we have
	\ben \E(H_M)-\E(H_K) = \int_{0}^{+\infty } \l(\int_0^t( L_K(s)-L_M(s))ds \r)dt.\een
	To prove \eref{eq:qdgr}, 
	it suffices to show that 
	\begin{lem}\label{lem:concl} We have \ben B_K(t) \geq B_M(t) ~~ \textrm{ for all }t\een
		where 
		\[B_K(t)=\int_0^t L_K(s)ds,~~B_M(t)=\int_0^t L_M(s)ds\]
		are the volume of $K$ and of $M$ below level $t$.
	\end{lem}  
	\begin{lem} \label{lem:zgthehgzfed}
		Let $K$ be an element of $\CCSF(F)$, and $m_K$ its maximal height.  
		For every $(t_1,s,t_2)$ such that $0\leq t_1\leq s\leq t_2\leq m_K$, we have
		\ben\label{eq:form1}
		L_K(s)^{\frac1{d-1}} \geq  \frac{t_2-s}{t_2-t_1}L_K(t_1)^{\frac1{d-1}} +\frac{s-t_1}{t_2-t_1}L_K(t_2)^{\frac1{d-1}},
		\een
		which is equivalent to say that $s\mapsto L_K(s)^{\frac1{d-1}} $ is concave. 
		For a unitary mountain $M$ (with maximal height $d$)  for $t\in[0,d]$,
		\[L_M(t) = \l(1-\frac{t}{d}\r)^{d-1}, ~~B_M(t)=1-\l(1-\frac{t}{d}\r)^{d},\]
		thus  $s\mapsto L_M(s)^{\frac1{d-1}}$ is linear, and in this case $K=M$, \eqref{eq:form1} becomes the equality case.
	\end{lem}

	\begin{proof}  
		Since $K$ is convex, any line intersecting both $\Layer_K(t_1)$ and $\Layer_K(t_2)$ intersects also $\Layer_K(s)$, and then 
		\ben \Layer_K(s)\supset \frac{t_2-s}{t_2-t_1}\Layer_K(t_1)+\frac{s-t_1}{t_2-t_1}\Layer_K(t_2),
		\een
		(since each element
		of the type  $ z:=\frac{t_2-s}{t_2-t_1} z_1+\frac{s-t_1}{t_2-t_1}z_2 $ for some $z_1\in \Layer_K(t_1)$ and $z_2\in\Layer_K(t_2) $ is in 
		$\Layer_K(s)$).
		By the Brunn-Minkowski inequality (see e.g. Schneider \cite[Section 6.1]{MR1216521})
		applied to the $(d-1)-$dimensional layers, we get 
		\ben\Vol_{d-1}^{\frac1{d-1}}\l(\Layer_K(s)\r)\geq \Vol_{d-1}^{\frac1{d-1}}\l(\frac{t_2-s}{t_2-t_1}\Layer_K(t_1)\r)+\Vol_{d-1}^{\frac1{d-1}}\l(\frac{s-t_1}{t_2-t_1}\Layer_K(t_2)\r)\een
		and since $\Vol_{d-1}(c\Layer_K(\tau))=c^{d-1}\Vol_{d-1}(\Layer_K(\tau))$, we get \eref{eq:form1}. The other statements are just exercises.
		
	\end{proof}

	To complete the proofs of the lower bound parts of  Theorems \ref{theo:Prism}  and  \ref{theo:CCSF}, it remains to prove Lemma \ref{lem:concl}. 
	
	\begin{proof}[Proof of  Lemma \ref{lem:concl}.]
		There is nothing to prove if $K$ is a unit mountain: we will  assume that $K$ is not a unit mountain, so that its maximum height satisfies $m_K<d$. 
		
		Consider the map $\Psi:\R^+\to \R$
		\[\Psi(t):=  L_K(t) ^{\frac{1}{d-1}}- L_M(t)^{\frac{1}{d-1}}.\]
		We claim that there exists $t^\star\in(0,m_K)$ such that $\Psi(t)$ is non-negative on $[0,t^\star]$, non-positive on $(t^\star,m_K)$. 
		
		The main argument will be the concavity of $t\mapsto L_K(t) ^{\frac{1}{d-1}}$ and of $t\mapsto - L_M(t) ^{\frac{1}{d-1}}$ (notice the sign): there is only a tiny complication that comes from the fact that the first one is concave on $[0,m_K]$ (and null thereafter), while the second is concave on the larger interval $[0,d]$ (and null thereafter).\par
		
		To prove the claim we need three facts:\\
		$\bullet$ By \eref{eq:form1}, $\Psi$ is concave on  $[0,m_K]$ (and negative on $(m_K,d)$, and null thereafter).\\
		$\bullet$ We have $\int_0^{m_K}L_K(t)dt=1>\int_0^{m_K}L_M(t)dt$, so that there must exist some $t$ such that $L_K(t)>L_M(t)$ in $(0,m_K)$, so that $\Psi$ must be positive somewhere on $(0,m_K)$.\\
		$\bullet$ In the extremities, we have $\Psi(0)=0$ and at this point there are two cases depending on the sign of $\Psi(m_K)$ : \\
		\indent-- either $\Psi(m_K)\geq0$, and since $\Psi(t)<0$ on $(m_K,d)$ the claim holds for $t^\star=m_K$,\\
		\indent-- or $\Psi(m_K)<0$, in which case, by concavity, there is a single cancellation time $t^\star$ on $(0,m_K)$, and $\Psi(t)$ is non-positive on $(t^\star,d]$.\\
		This proves the claim.
		
		Since the map $x\mapsto x^{1/(d-1)}$ is non-decreasing (on $\R^+$), the claim is also valid for the map  $\Phi(t):= L_K(t)-L_M(t)$.
		
		Therefore the map $t\mapsto B_K(t)-B_M(t)$, whose derivative is $\Phi(t)$, is non-decreasing on $[0,t^\star]$ and non-increasing on $[t^\star,d]$. Since $B_K(0)-B_M(0)=B_K(d)-B_M(d)=0$, these variations properties imply that  $t\mapsto B_K(t)-B_M(t)$ is non-negative on $[0,d]$ (and null after that).
	\end{proof}

	\subsubsection{Proof of Theorem \ref{theo:CCSF} : upper bound}
	
	The main idea is to construct a wide trapezoid $T_h$ with a very small height, so that uniform random point inside have small heights. 
	The first base of the trapezoid is $B_1:=F\times\{0\}$ and the second one is  $B_2:=(c_h F)\times \{h\}$ for $h>0$, $c_h$ a real positive number (depending on $h$), and the trapezoid $T_h$ is the convex hull of the two bases (and has floor $F$). 
	Since $\Vol_{d-1}(B_2)=c_h^{d-1}$, and since we want $\Vol_d(T_h)=1$, we set $c_h = (2/h-1)^{1/(d-1)}$. 
	When $h$ goes to zero, $c_h$ goes to $+\infty$. Take $\alpha>0$.
	For $U=(Z,H)$ uniformly taken in $T_h$, $\E(H)\leq h$ and then $Q_{T_h}(2)= 1 - 2\E(H)/d\geq 1-2h/d\geq 1- \alpha$, for $h>0$ small enough. 
	
	It remains to show that $Q_K(2)=1$ is impossible: since $K$ has a non-empty interior, for $U=(Z,H)$ uniform in $K$, $H$ is almost surely positive, and then $\E(H)>0$.

	\section{2D - considerations}
	
	In the 2D case, we will make an extensive use of top functions, which are simple concave functions. We focus on convex sets having floor $F:=[0,1]\times\{0\}$, belonging to $\SubPrism(F),$ which  prosaically means that they are contained in the rectangle 
	\[{\sf Rect}:=[0,1]\times[0,2].\]
	
	We denote by  $\Conca$ the set of concave functions $G:[0,1]\to \R^+$ such that  $\int_0^1  G(x)\,dx=1$. The set $\Conca$ and $\SubPrism(F)$ can be identified by associating with each $K$ its top function and reciprocally, to a function $G\in\Conca$, the set below its graph 
	\[D_G:=\{ (x,y)~: x\in[0,1], y\leq G(x)\}\cap {\sf Rect}\] (see if needed the beginning of Section \ref{sec:egther}).
	Since it is more convenient to work with functions in the sequel, we will use the notation
	\[Q^{G}_n:= Q_{D_G}(n).\]

	\subsection{General decomposition formula for $Q_K(n)$}
	\label{GDF}
	
	The main message of this section is that there exists a formula that allows to compute $Q^G_n$: the formula is recursive and uses some elements $Q^{H}_n$ where $H$ are defined using a part of the graph of $G$. Let us introduce the material devoted to the expression of this decomposition.
	
	Take $t\in[0,1]$ and consider the two segments $S_1:=[ (0,0),(t,G(t))]$ and $S_2:=[(t,G(t)),(1,0)]$. Consider $L(t)$, the points that are in $D_G$ and above the segment $S_1$; similarly, let $R(t)$ those that are in $D_G$ and above $S_2$. Formally, 
	\ben
	L(t) &=&\l\{(x,y)~: 0\leq x \leq t, y\geq \frac{x}{t} G(t) \r\} \cap D_G,\\
	R(t) &=&\l\{(x,y)~: t\leq x \leq 1, y\geq G(t)+\frac{x-t}{1-t} (G(1)-G(t)) \r\} \cap D_G. 
	\een
	In what follows, we will sometimes write $L^G(t)$ and $R^G(t)$ to clarify the dependence of these sets on $G$. For any $E\subset\R^2$, we will also use the notation $\abso{E}$ standing for the Lebesgue measure of the set $E$ ; hence 
	\[|L(t)|+|R(t)|=1-G(t)/2.\]
	\begin{figure}[htbp]
		\centerline{\includegraphics{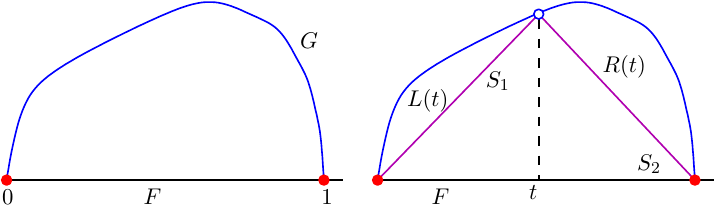}}
		\caption{Introduction of the domains $L(t)$ and $R(t)$}
	\end{figure}
	
	In the sequel, we perform a decomposition at any $t\in[0,1]$ along $L(t)$ and $R(t)$ and deal with the uniform distribution in these domains which have respective floors $S_1$ and $S_2$. 
	
	We would rather work with the fixed floor $F=[0,1]\times\{0\}$ and functions in $\Conca$, therefore we map these domains at a suitable position using affine maps (that preserve both  convexity and the uniform distribution). Note that for $t\in(0,1)$ one can find:\\
	$\bullet$ a unique affine map $\Phi^{L,t}$ with positive determinant, that sends $S_1$ on $F$ and $L(t)$ inside ${\sf Rect}$, preserves vertical lines, and such that $\big|\Phi^{L,t}(L(t))\big|=1$; let $NL(t)$ be the top function of $\Phi^{L,t}(L(t))$ (the prefix "$N$" stands for normalized)  \\ 
	$\bullet$ a unique affine map $\Phi^{R,t}$ with positive determinant, sending $S_2$ on $F$ and $R(t)$  in ${\sf Rect}$, preserving vertical lines, and such that $\big|\Phi^{R,t}(R(t))\big|=1$. Let $NR(t)$ be the top function of $\Phi^{R,t}(R(t))$.

	The sequence   $(Q^G_n,n\geq 0)$ obeys the following decomposition
	\begin{theo}\label{theo:peigne}For any function $G \in \Conca$, any $n\geq 0$,  
		\ben\label{eq:grhet}
		Q^G_n =  \sum_{m=0}^{n-1} \binom{n-1}m \int_0^1 G(t) Q^{NL(t)}_m Q^{NR(t)}_{n-1-m} |L(t)|^{m}  |R(t)|^{n-1-m} dt.
		\een 
	\end{theo}
	\begin{rem} Theorem \ref{theo:peigne} allows to express $Q^G(n)$ in terms of $Q^{NR(t)}_m$ and $Q^{NL(t)}_{n-1-m}$ where $NR(t)$ and $NL(t)$ are relatively simple functions of the initial function $G$. In order to derive a closed formula for $Q^G_n$, it is morally necessary that the class of domains $NL(t)$ and $NR(t)$ appearing in the decomposition remains under control. This is the case if the class of functions is stable by composition by affine functions (as constant, linear, quadratic functions are). 
	\end{rem}
	\begin{proof}The proof is inspired by the comb formula in \cite{marckert:hal-02913348}.
		
		Consider the family of maps $G_\lambda:x \mapsto (1-\lambda)G(x)$. For $\lambda=0$ the graph of $G_{\lambda}$ is the same as that of $G$. As $\lambda$ moves from 0 to $1$, the family of domains $D^{G_\lambda}$ decreases for the inclusion partial order, and eventually, at time $\lambda=1$, becomes flat.
		\begin{figure}[htbp]
			\centerline{\includegraphics{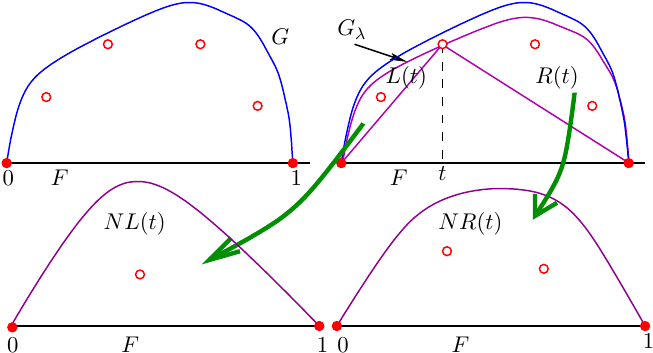}}
			\caption{\label{dif:decGlambda} Families of map $G^{\lambda}$ (decreasing as $\lambda$ grows). Representation of the images $NL(t)$ and $NR(t)$ of $L(t)$ and $R(t)$ by the affine maps that transport the "floor segments" of both domains onto $F$, that send vertical lines onto vertical lines, and that vertically rescale the domain so that $|NR(t)|=|NL(t)|=1$. }
		\end{figure}
		Let $U_1,\cdots,U_n$ be iid uniform in $D^G$. Let $\lambda^{\star}=\sup\{\lambda>0: U_1,\cdots,U_n \in D^{G_\lambda} \}$. Almost surely $\lambda^\star\in(0,1)$, and for $\lambda^\star=\lambda$ there is a single (uniform index) $i^\star\in\cro{1,n}$ such that $U_{i^\star}\in\partial D^{G_\lambda}$.
		Now, the set of points $U_1,\cdots,U_n$ together with $(0,0)$ and $(1,0)$ are in convex position if and only if none of them lies in the interior of the triangle with vertices $(0,0)$, $(1,0)$, $U_{i\star}$, and if the points of $\cro{1,n}\setminus\{i^\star\}$ can be partitioned in two parts $P_1, P_2$ such that all the points with indices in $P_1$ lie in $L^{G_\lambda}(t)$ and are in convex position together with the extremities of $S_1$ (computed in $D^{G_\lambda})$, and all those of $P_2$ lie in $R^{G_\lambda}(t)$ and are in convex position together with the extremities of $S_2$ (computed in $D^{G_\lambda}$).

		We claim that, for $U_i=(X_i,Y_i)$
		\ben\label{eq:qrjdgr}
		Q_n^G = \sum_{k=1}^n\sum_{m=0}^{n-1} \int_{0}^1 \int_{0}^1 \binom{n-1}m |L(t)|^m|R(t)|^{n-m-1} Q^{NL(t)}_m Q^{NR(t)}_{n-m-1}\P(i^\star=k, X_{i^\star} \in dt)   d\P_{\lambda^\star}(\lambda)
		\een
		where  $\P_{\lambda^\star}$ stands for the law of $\lambda^\star$.
		This formula is valid because:\\
		\indent-- when one finds $U_{i^\star}$ in $\partial G_{\lambda^\star}$, the other points are uniformly distributed in $D^{G_{\lambda^\star}}$. We search for the probability that these points are in convex position (together with $(0,0)$ and $(1,0)$ and $U_{i^\star}$), and avoid the triangle $(0,0), (1,0)$, $U_{i\star}$. By applying the affine map $(x,y)\to (x, y/(1-\lambda))$ that sends $D^{G_\lambda}$ over $D^{G}$, we are back to the  original problem with $n-1$ points that must fall partly in $L(t)$, and partly in $R(t)$: this means that the initial function $G$ still plays a role in the definition of the left and right domains $L(t)$ and $R(t)$ even if we had contracted $D^G$ into $D^{G_\lambda}$.\\
		Now, we claim that 
		\ben\label{eghtyjthfd}\P(i^\star=k, X_{i^\star} \in dt)  =(1/n) G(t) dt,\een
		where here and everywhere, $\P(X\in dt)=f(t) dt$ means that law of the random variable $X$ has density $f$ with respect to the Lebesgue measure.
		Indeed, in order to sample $U_i=(X_i,Y_i)$ uniformly at random in $D^G$, one can proceed as follows: sample $X_i$ according to $G$ and then set $Y_i=G(X_i)V_i$ where $V_i$ is uniform on $[0,1]$ (all the $V_i$ are independent, and independent of the $X_i$). In this representation, $i^\star$ is the index $i$ of the largest $V_i$, and the remaining $X_i$ are still independent and have distribution $G$.
		Injecting \eref{eghtyjthfd} in Formula \eref{eq:qrjdgr}, allows to see that the sum over $k$ can be simplified, and the integral over $\lambda$ can be factorized, and gives 1.
		
	\end{proof} 
	
	\subsection{Exact computation for some particular shapes}
	
	\subsubsection{The bi-pointed triangular case} 
	Consider the triangle $\triangle$ defined by the concave map of area 1, $G(t)=2t$ (on $[0,1]$), with floor $F=[0,1]\times\{0\}$. 
	Let 
	\[t_n=Q_\triangle(n)\] be the corresponding bi-pointed probability.
	Let us recover the following theorem due to Bárány, Röte, Steiger and Zhang \cite{baranybipointed} (see also Gusakova and Th\"ale \cite{GT2021} for recent results on the number of points on the convex hull of $\{U_1,\cdots,U_n\} \cup\{(0,1),(0,0)\}$ where the $U_i$ are iid uniform random points $n$ in the triangle).
	\begin{theo}\label{theo:tri} For all $n\geq 1$, we have
		\ben\label{eq:tn} t_n = \frac{2^n}{n!(n+1)!}\een
	\end{theo}
	The eight first terms of this sequence (starting with $t_0$) are
	\ben
	1,1,{\frac{1}{3}},{\frac{1}{18}},{\frac{1}{180}},{\frac{1}{2700}},{
		\frac{1}{56700}},{\frac{1}{1587600}},{\frac{1}{57153600}}
	\een
	Recall that for all $a,b>0$,
	\ben\label{eq:beta}B(a,b):=\int_{0}^1 t^{a-1}(1-t)^{b-1}dt= \Gamma(a) \Gamma(b) / \Gamma(a+b).\een 
	
	\begin{proof} 
		
		For each $t$,  $L(t)$ is flat and $R(t)$ is triangular and has area 
		${|R(t)|}= 1-t$. We then have
		\[t_n = \int_{0}^1 (2t) (1-t)^{n-1} t_{n-1} \binom{n-1}{n-1} dt=2B(2,n) t_{n-1}\]
		and using \eref{eq:beta}, we conclude (since $B(2,n)=1/(n(n+1))$). 
	\end{proof}

	\subsubsection{The bi-pointed square case}
	
	In a square $S =[0,1]^2$ with floor $F=[0,1]\times\{0\}$, set $q_n=Q_S(n)$. 
	\begin{pro}\label{pro:q} For all $n\geq 0,$ we have
		\ben\label{eq:qn} q_n =\frac1{n!(n+1)!}{2n \choose n}.\een
	\end{pro}
	The eight first terms are 
	\[1,1,{\frac{1}{2}},{\frac{5}{36}},{\frac{7}{288}},{\frac{7}{2400}},{
		\frac{11}{43200}},{\frac{143}{8467200}},{\frac{143}{162570240}}\]
	
	\begin{proof}
		We have $G(t)=1$ for all $t\in[0,1]$, $R(t)$ is a triangle with area $t/2$, $L(t)$ a triangle with area $(1-t)/2$. Still using  \eref{eq:beta}, we have:
		\ben q_n &=& \sum_{k=0}^{n-1} \binom{n-1}k t_k t_{n-1-k} \int_{0}^1  (t/2)^k ((1-t)/2)^{(n-k-1)}   dt\\
		&=& \sum_{k=0}^{n-1} \binom{n-1}k \frac{ 1}{(k+1)!(n-k)!} \frac{1}{n!}=\frac{1}{n!(n+1)!}\sum_{k=0}^{n-1} \binom{n-1}k \binom{n+1}{n-k};
		\een
		to conclude, observe that this last sum is $\binom{2n}n$ because one can partition the set of subsets of size $n$ of $\{1,\cdots, 2n\}$, with respect to their numbers $k$ of elements belonging to  $\{0,\cdots,n-1\}$ (and then, $k$ ranges from 0 to $n-1$).
	\end{proof}
	
	\subsubsection{The bi-pointed parabola case}
	
	Consider the parabolic function $G:x\mapsto 6 x(1-x)$; its area is 1, and of course, $G$ is in $\Conca$.
	Let $p_n=Q_n^G$ be the probability that $n$ iid  point taken uniformly in $D$, together with $(0,0)$ and $(1,0)$ are in convex position. 
	\begin{theo}\label{theo:para} For all $n\geq 0$ we have
		\ben\label{eq:pn} p_n = \frac{12^{n+1}}{6 (2 n+2)!}.\een 
	\end{theo}
	\begin{proof} In this case,  
		\[|L(t)|= t^3\textrm{ and } |R(t)|=(1-t)^3.\]
		Moreover, $L(t)$ (resp. $R(t)$) can be sent on $D^G$ through  $\Phi^{L,t}$ (resp.  $\Phi^{R,t}$), so that $NL(t)=NR(t)=G$.
		We then have by Theorem \ref{theo:peigne}, $p_0=1$, and for $n\geq 1$,
		\ben p_n = \sum_{k=0}^{n-1} \binom{n-1}k p_k p_{n-1-k} \int_{0}^1 G(t) t^{3k} (1-t)^{3 (n-k-1)}dt.
		\een
		Replacing $G(t)$ by $6t (1-t)$ (inside this integral) and using  \eref{eq:beta},
		we get 
		\ben  p_n = \sum_{k=0}^{n-1} \binom{n-1}k p_k p_{n-1-k} \times 6 \frac{(1+3k)! (1+3(n-k-1)!) }{(3n)!}.\een
		To prove that the sequence $(p_n)$ given in \eref{eq:pn} satisfies this formula, it suffices to assume that it does so for all $k<n$ (for some $n\geq 1$), and to prove that it does so for $n$ as well.
		To do this, replace in the rhs $p_k$ by $12^{k+1}/(6(2k+2)!)$.
		We obtain 
		\begin{align*}
			p_n &= \frac{12^{n+1}}{6}\frac{(n-1)!}{(3n)!}\cdot\sum_{k=0}^{n-1} \frac{(1+3k)! (1+3(n-k-1)!) }{k!(n-1-k)!(2k+2)!(2n-2k)!},
		\end{align*}
		so that it suffices to prove that 
		$$\sum_{k=0}^{n-1} \frac{(1+3k)! (1+3(n-k-1)!) }{k!(n-1-k)!(2k+2)!(2n-2k)!}=\frac{(3n)!}{(n-1)!(2n+2)!}$$
		to conclude. Rearranging the terms, this formula is equivalent to
		\ben\label{eq:fgr} 
		\binom{3(n-1)+4}{n-1}\frac{4}{(3n+1)}= \sum_{k=0}^{n-1} \binom{3k+2}{k} \frac{2}{3k+2}  \binom{ 3(n-k-1)+2}{(n-k-1) } \frac{2}{3(n-k-1)+2}\een
		or in other words to
		\ben \label{eq:htime} {\sf Q}[3(n-1),4] = \sum_{k=0}^{n-1} {\sf Q}[3k,2]\; {\sf Q}[3(n-k-1),2]\een
		where
		\[{\sf Q}[ 3N,F ] :=\binom{3N+F}{N} \frac{F}{3N+F}.\]
		
		The proof of Theorem \ref{theo:para} will be completed when we will have proven Formula \eref{eq:htime}, which is a well-known consequence of the ``cycle lemma'' (for $(S_1,S_2)=(2,2)$), recalled below for the convenience of the reader.
	\end{proof}
	\begin{lem} Take $N>0$, $S>0$. 
		The number of paths in $\Z$ with step $+2$ or $-1$, starting at 0, of length $3N+S$, that hits $-S$ for the first time at time $3N+S$ is counted by $\binom{3N+S}{N} \frac{S}{3N+S} = {\sf Q}[3N,S]$. Hence , for any $S>0$, $S_1>0$, $S_2>0$ such that $S=S_1+S_2$,
		\[{\sf Q}[3N,S] = \sum_{k=0} ^{N}{\sf Q}[3k,S_1]{\sf Q}[3(N-k),S_2].\] 
	\end{lem}
	\begin{proof}
		This Lemma is a consequence of a classical result in combinatorics which appears in the literature under many names: the cycle lemma,   Dvoretzy-Motzkin's lemma \cite{DM},  Otter formula, Kemperman' s formula, and so on (see Pitman \cite[Section 6.1]{MR2245368} for historical notes and discussion).
		
		It possesses a probabilistic and a combinatorial formulation, as follows:
		Consider $AI$ a set of ``allowed increments'' included in $\{ -1, 0, 1,2,...\}$ containing at least $-1$ and a positive element.  Denote by ${\sf Paths}[n, f] = \{ w[n]=(w_0,\cdots,w_{n}):~:~ w_i-w_{i-1}\in AI, w_0=0,w_n= -f \}$ the set of paths of length $n$ whose increments are allowed, starting at zero, ending at $-f$, for some $f>0$. For a path $w[n]$, denote by $\tau_f(w[n])$ the first index $k$ such that $w_k=-f$ (the hitting time of $-f$). 
		Let 
		\[{\sf Path}^{\star}[n,f]=\{ w[n] \in {\sf Paths}[n, f], \tau_{f}(w[n])=n\}\] that is, the elements of ${\sf Paths}[n, f]$ hitting $-f$ at the end, at time $n$.  The cycle lemma asserts that if $n\geq 1$ and $ f>0$, 
		\[|{\sf Path}^{\star}[n,f]| = (f/n)\,\, |{\sf Path}[n,f]|.\] 
		The probabilistic version (proved exactly the same way) is valid when the increments are iid, with a distribution having supports in $AI$. We then sample a random walk $W=(W_0,\cdots,W_n)$ with these increments, and denoting by $\tau_f$ the hitting time by the random walk of $-f$, the formula goes as $\P(\tau_{f}(W[n])=n)=(f/n) \P(W_n=-f)$. 
		
		Here, we take $AI=\{2,-1\}$.  In this case, observe that if $w$ is in ${\sf Paths}[n, f]$ for a fixed $f=S$, then, if $n_2$ and $n_{-1}$ are respectively its numbers of $+2$ and $-1$ steps, then
		\[ (n,-S) = (n_2 + n_{-1}, 2n_2 -n_{-1}),~~\imp (n_2,n_{-1}) =( \frac{n-S}{3},  \frac{2n +S}{3}).\]
		In particular, $n$ has the form $3N+S$ and $(n_2,n_{-1})=(N,2N+S)$ (for $N=(n-S)/3$).
		It follows that ${\sf Path}[3N+S,S]=\binom{3N+S}{N}$, and by the cyclic lemma   ${\sf Path}^{\star}[3N+S,S]=(S/(3N+S))\binom{3N+S}{N}$, which proves then the first statement of the Lemma.
		
		The second statement is easier: we will associate bijectively, to each element $w$ in ${\sf Path}^{\star}[3N+S,S]$, a pair of paths $(w^{(1)},w^{(2)}) $, where $w^{(1)}$ is in ${\sf Path}^{\star}[3k+S_1,S_1]$ and $w^{(2)}$ is in  ${\sf Path}^{\star}[3(N-k)+S_2,S_2]$, for some $k$. This association is a bijection from ${\sf Path}^{\star}[3N+S,S]$ onto $\cup_{k=0}^N{\sf Path}^{\star}[3k+S_1,S_1] \times {\sf Path}^{\star}[3(N-k)+S_2,S_2]$.
		
		Take a path $w$ in ${\sf Path}^{\star}[3N+S,S]$. It hits level $-S$ at time $3N+S$ for the first time. Since the only negative steps in $w$ are $-1$, this path hits level $-S_1$ for the first time  at some time $\tau_{S_1}$ of the form $3k+S_1$, and the increments of the rest of the paths correspond to a path with length $3N+S-(3k+S_1)=3(N-k)+S_2$ hitting $-S_2$ at the end. Technically:
		\be
		w^{(1)}_j &=& w_j , \textrm{ for } 0 \leq j \leq \tau_{S_1},\\
		w^{(2)}_j &=& w_{\tau_{S_1}+j}-w_{\tau_{S_1}},  0 \leq j \leq 3N+S- \tau_{S_1}=3(N-\tau_{S_1})+S_2.
		\ee
		The inverse of this bijection consists just to take a pair of paths  $(w^{(1)},w^{(2)})$  in ${\sf Path}^{\star}[3k+S_1,S_1] \times {\sf Path}^{\star}[3(N-k)+S_2,S_2]$ and to concatenate their increments in order to construct a path in ${\sf Path}^{\star}[3N+S,S]$.
	\end{proof}

	\begin{rem}B\'ar\'any's analogues in the $2D$-floor case?
		It is a well-known result, due to B\'ar\'any \cite{barany1}, that in any convex set $K$ with non-empty interior, the convex hull of $n$ iid uniform points taken in $K$, conditioned to be in convex position, converges for the Hausdorff distance $d_H$ in probability to a (unique) deterministic convex set $\Dom{K}$ inside $K$. The set $\Dom{K}$ can be characterized  both  as the aforementioned limit, or as the only set among all convex subsets of $K$ realizing the maximum of a quantity called the affine perimeter $\AP$ (we refer to \cite{barany1} for one of the definitions of the affine perimeter). We may conjecture that this property is still true in the $2D$-floor case, with the following adjustment: the limiting shape would be the convex set $C$ included in $K$ and containing the floor, that maximises the affine perimeter. Even if the arguments of B\'ar\'any seem to extend to the present models, a complete formal proof, would be needed.
	\end{rem}

	\subsection{Optimization of \texorpdfstring{$Q_K(2)$}{Q_K(2)} in 2D: a second simple proof}\label{sec:secondproof}

	In 2D, since we want to prove $Q_2^G \geq t_2$ where $t_2$ corresponds to $Q_2^H$ with $H:t\mapsto 2t$, it suffices to prove that
	\[Q_2^G= \int_0^1 G(t) (1 - G(t)/2) dt \geq \int_0^1 2t (1 - 2t/2) dt=t_2~~~\textrm{ for }G\in\Conca\]
	which is then equivalent to
	\begin{align*}
		\int_0^1 G^2(t) dt \leq \int_0^1 (2t)^2dt=\frac{4}{3}.
	\end{align*}
	
	The idea of this second argument relies on a ``representation formula'': every function $G\in \Conca$ can be written as a mixture of mountains -- that is a weighted sums of mountains -- that is, it can be written as
	\ben\label{eq:herytu2} G= \int_0^1   M_s d\mu(s)\een where $\mu$ is a positive Borelian measure with total mass 1 on $[0,1]$. Such a $\mu$ is a probability distribution on $[0,1]$, which means that for $S$ a random variable with law $\mu$,
	\[G(x)=\E( M_S(x)),~~~\textrm{ for all } x\in[0,1].\]

	Assume for a moment that such a representation exists and observe that for any mountain $M_s$, 
	\[\int_0^1 M_s^2(x) dx=\int_0^1 (2t)^2\, dt=4/3.\]
	Now, by Fubini
	\[\int_0^1 G(x)^2dx=\int_{0}^1 \E(M_S(x))^2 dx\leq \int_0^1 \E(M_S(x)^2) dx=\E\l(\int_0^1 M_S(x)^2 dx\r)\leq \int_0^1 \frac 43 dx = \frac 43.\] 
	It remains to establish the representation formula, and we will do so by proving it first in the case where $G$ is piecewise linear, with a finite number of pieces.
	In this case  $G'(a^+)-G'(a^-)$ is zero except at $a\in I $  for a finite set $I$ , corresponding to the abscissas of the vertices of the polygonal graph of $G$, corresponding to discontinuities in the derivative of $G$.
	Consider a sum of mountains with same set of vertices: $H:=\sum_{s\in I \cup\{0,1\}} p_sM_s$ where the $p_s$ are non negative (and positive if $s\in S$). Since for the mountains $M_s$ and $s\in(0,1)$,
	\[M'_{s}(s^+)-M'_s(s^-) =-2/(1-s) - 2/s\] it is possible to fix the $p_s$ so that 
	\[G=\sum_{s\in I\cup \{0,1\}} p_s M_s\]
	by taking $p_s =  \frac{ G'(s^+)-G'(s^-)}{-2/(1-s) - 2/s}$ for $s\in I$, and $p_0=G(0)/2$, $p_1=G(1)/2$.
	The fact that $\sum_{s\in I\cup\{0,1\}} p_s=1$ is just a consequence of the fact that $\int G=\int H=1$, and $\int M_s=1$ for all $s$.
	This provides a representation formula for these types of functions $G$, which are dense in \Conca, for the uniform topology. For a general $G$ in $\Conca$, take a sequence $(G_n)$ such that $G_n$ is piecewise linear (with a finite number of pieces), and $G_n\to G$ uniformly. The corresponding sequence of measure $(\mu_n)$ possesses an accumulation point $\mu$ in the set of probability measures on $[0,1]$ (which is compact) and then there exists a converging subsequence $(\mu_{n_k})$, that converges to a probability measure $\mu$ on $[0,1]$. Let $Z$ be $\mu$ distributed. Now, for all $x\in[0,1]$,
	\[G(x)=\lim G_n(x) =\lim_k G_{n_k}(x) = \lim_k \E\l(M_{S_{n_k}}(x)\r).\]
	Since the map $s\mapsto M_s(x)$, is bounded and continuous on $[0,1]$, we get from $S_{n_k}\dd Z$,
	\[G(x)= \E(S_X(x))=\int_0^1 M_s(x) d\mu(s),\]
	which is the wanted representation formula.
	
	\section{Convex sets in 3D with a floor: new results}
	
	In 3D, even for the seemingly simplest case as the cube, the tetrahedron with a 2 dimensional floor, we have no formula for $Q_K(n)$ (but $Q_K(2)$ can be computed in these cases). To get an idea on the behavior of $Q_K(n)$,  we try to give insights and bounds on this quantity when $K$ is either a prism (that we will denote $P$ for short) or any mountain $M$.

	\subsection{A universal lower bound of \texorpdfstring{$Q_K(n)$}{Q_K(n)} for mountains in 3D}
	
	\begin{theo}\label{theo:dfreghrtu} 
		Let $d=3$. Take any floor $F$ (compact convex subset of $\R^2\times\{0\}$ with area 1), any  unit mountain $M$ with floor $F$. For all $n\geq 0$,
		\[Q_M(n)\geq Y_n:=\frac{2^n}{n!}\prod_{j=1}^n \frac{1}{3j-1}.\]
		The first terms of the sequence $(Y_n)$, for $n\geq 0$, are the following:
		\[1,1, \frac{1}{5}, \frac{1}{60}, \frac{1}{1320}, \frac{1}{46200}, \frac{1}{2356200}, \frac{1}{164934000}, \frac{1}{15173928000}.\]
	\end{theo}
	We will assume all along the proof that the center of mass of $F$ is at (0,0,0).

	The idea is as follows: put $U_1,\cdots,U_n$ in $M$, and consider the respective heights of the layers $h_1,\cdots,h_n$ containing these points (note that these heights cannot exceed $h=d=3$). The $\Layer_h$ is a copy of  $\Layer_0$, rescaled, and placed at level $h$. We have  \be\Layer_M(h)&=& (1-h/3) F + (0,0,h)\\
	\Layer_P(h)&=&   F + (0,0,h).
	\ee
	The layers are copies of $F$. If a point $z$ (random or not) falls in a layer, we may associate a factor $a(z)$, that we will call $F$-radius, which measures the smallest copy of $F$ (still with center of mass $0$) which would contain it: for $z\in \Layer(h)$:
	\[a(z)=\min\{ a : a \geq 0, z \in (aF) + (0,0,h)\}.\]
	The knowledge of $(a(z),h)$ is not sufficient to recover $z$, since given these numbers, we only know that $z$ belongs to the boundary of $(a(z)F) + (0,0,h)$.

	\begin{figure}[htbp]
		\centerline{
			\includegraphics[height=7.2cm]{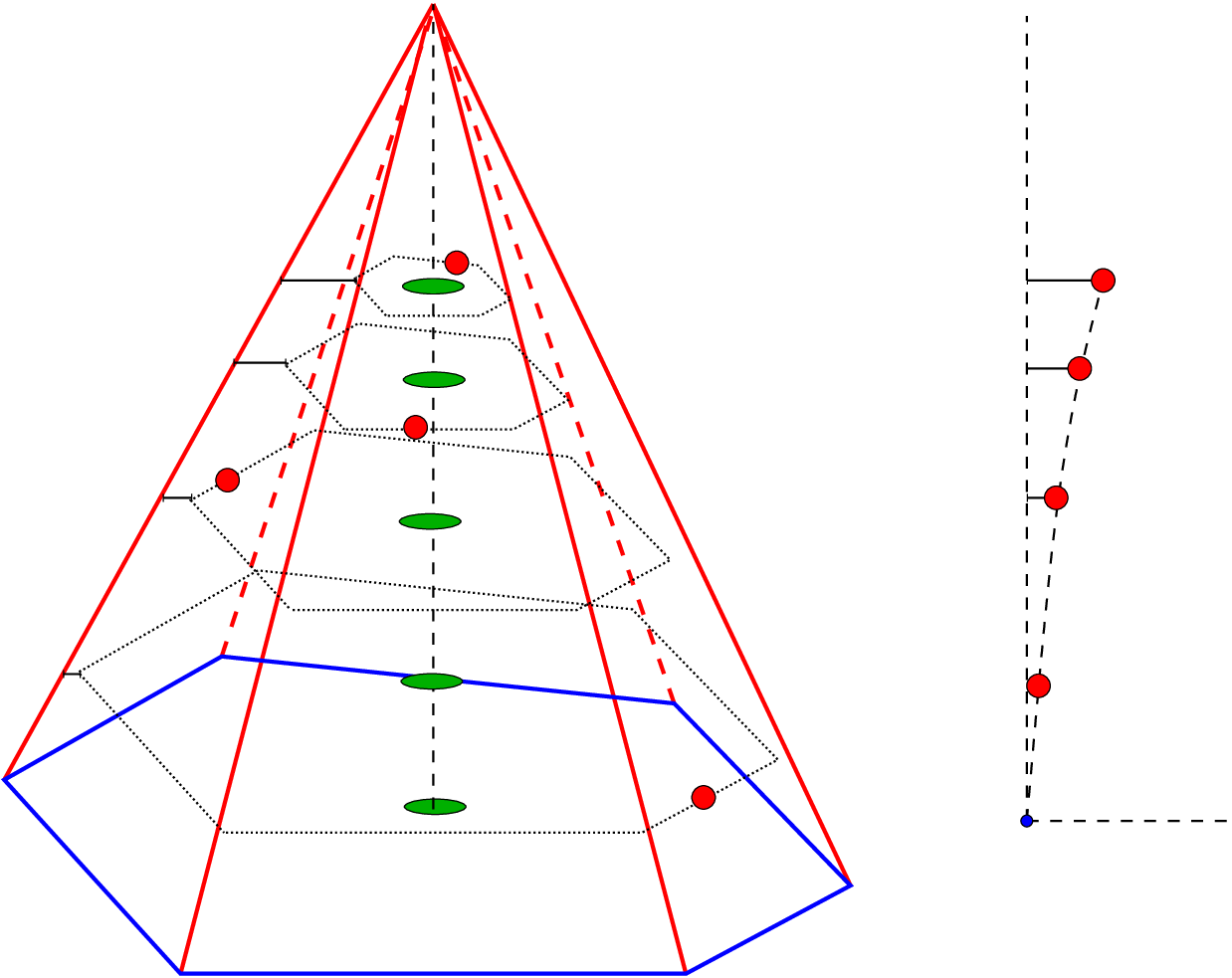}}
		\caption{An illustration of Lemma \ref{lem:boundsg}: on the left, the points $(z_i)$ are in convex position above the polygonal floor. Each horizontal slice of the mountain containing a point is a dilatation $F'$ of $F$ (with center of mass above that of $F$), and itself contains a second dilatation $F''$ of $F$ (with center of mass above that of $F$) having this point on its boundary (these are the represented by gray lines). The scale factor needed to go from $F$ to $F''$ is $a$.  On the right figure are represented the 2D points $(a_i,h_i)$ of the successive heights and factors of the points, from the bottom to the top. If the 2D points $(a_i,h_i)$ are in convex position (together with $(1,0)$), then the $(z_i)$, together with the floor, are in convex position.}
		\label{fig:lemboundsg}
	\end{figure}
	\vspace{0.3cm}
	We now state and prove two structural lemmas that will allow us to prove Theorem \ref{theo:dfreghrtu} .
	In the sequel, for a point $z=(z_1,z_2,z_3)\in\R^3$, we let $\pi_3:z\mapsto z_3$.
	{\begin{lem}\label{lem:boundsg}For some valid floor $F$, consider a sequence of points $z_1,\cdots,z_n$ in $F\times [0,3]$ ordered by increasing third coordinates (\ie such that  $0\leq \pi_3(z_1)\leq \cdots \leq \pi_3(z_n)$). If the points 
			\[\big(1,0\big), \big(a(z_1),\pi_3(z_1)\big), \cdots, \big(a(z_n),\pi_3(z_n)\big)\] are in convex position in the 2D-rectangle $[0,1]\times[0,3]$, then the $(z_i)_{1\leq i\leq n}$ together with the floor $F$, are in convex position in $\R^3$.
		\end{lem}
		
		\begin{proof}   The point $(1,0)$ can somehow be considered as a point $z_0$ in $\partial F$ (it encodes the same information as $(a(z_0),\pi_3(z_0))$).
			
			Take a sequence $(a_i,h_i)_{0\leq i \leq n}$ with $(a_0,h_0)=(1,0)$, the $h_i$ are non decreasing and bounded by $3$, the $a_i$ non increasing, non negative, and the $(a_i,h_i)$ form a convex chain in $[0,1]\times [0,3]$. Now, consider the union of ``1-D'' curves
			\[ S_n:=\cup_{i=0}^n  (0,0,h_i)+a_i \partial F.\]
			We claim that $S_n$ is in convex position with $F$ (meaning that $F$, as well as each curve $(0,0,h_i)+a_i \partial F$ are on the boundary of $\CH(F \cup S_n)$ (the reason is that it is easy to construct a supporting plane at each point of $S_n$ using a supporting line at the corresponding point in $F$).
			
			Now, the conclusion follows because $(z_i)_{0\leq i\leq n}$ is a subset of $S_n$. \end{proof}}
	
	Pick $U_1,\cdots,U_n$ iid and uniform in $M$ or in $P$, and let us get to \eqref{theo:dfreghrtu}
	
	\paragraph{Comparison with "2D-problems".} The idea is now to identify sets of 3D points in convex position such that the Lebesgue measure of their set of 2D projections is computable.
	
	Consider the densities 
	\ben g_1(x,y)&=&2x \textrm{ on }[0,1]\times[0,1],\\ g_2(x,y)&=&2x \1_{x\leq (1-y/3)} \textrm{ on }[0,1]\times[0,3].\een
	\bls Denote by $\beta_1(n)$ the probability that $(1,0),(X_1,Y_1),\cdots,(X_n,Y_n)$ are in convex position in  $[0,1]\times[0,1]$, where the $(X_i,Y_i)$ are iid and have density $g_1$.\\
	\bls And similarly let $\beta_2(n)$ be the probability that $(1,0),(X_1,Y_1),\cdots,(X_n,Y_n)$ are in convex position in  $[0,1]\times[0,3]$, where the $(X_i,Y_i)$ are iid and have density $g_2$.
	
	\begin{lem}\label{lem:relationqes} For all $n\geq 0$,
		\be
		Q_P(n) & \geq & \beta_1(n),\\
		Q_M(n) & \geq & \beta_2(n).
		\ee
	\end{lem}

	\begin{proof}[Proof of Lemma \ref{lem:relationqes}] We now use Lemma \ref{lem:boundsg}, and seek to describe the joint laws of $(a(U_i),\pi_3(U_i))$ for uniform point $U_i$, in the two cases of interest.\par
		\bls In the prism case, $\pi_3(U_i)$ is uniform on $[0,1]$ and independent of $a(U_i)$. Now, conditional on $U\in\Layer_P(h)$, $U$ is uniform on this layer, which is isomorphic to $F$. To compute the law of $a(U)$, we may  work directly on $F$: for $a\in[0,1]$, $\P(a(U)\leq a)=\P(U\in aF)=a^2$, and then the density of $(a(U),\pi_3(U))$ is indeed $g_1$.
		
		\bls In the mountain case, $\pi_3(U_i)$ has density $(1-y/3)^2$ on $[0,3]$, since this is the height density, already discussed several times before. Now, given that $U$ is in $\Layer_M(y)$, it is again uniform in this layer, which is isomorphic to $(1-y/3)F$, which has $F$-radius $(1-y/3)$, and then $a(U)$ will be smaller than this value. We have $\P(a(U)\leq a ~|~ U\in \Layer_M(y))= a^2/((1-y/3)^2)\1_{a\leq 1-y/3}$ so that finally,  $\P((a(U),U)\in (dx,dy))= \frac{2a}{(1-y/3)^2} (1-y/3)^2 \1_{x\leq 1-y/3} dxdy=g_2(x,y)dxdy$ as announced.
	\end{proof}
	
	\begin{proof}[Proof of Theorem \ref{theo:dfreghrtu} ] 
		\begin{figure}[htbp]
			\centerline{\includegraphics{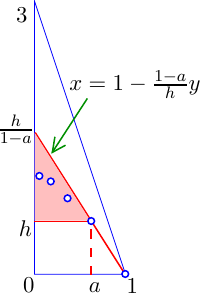}}
			\caption{\label{fig:lafig} On this representation, one sees that when the point $(a,h)$ with the smallest height $h$ is found, the other points must be in the colored zone to be in convex position.}
		\end{figure}
		We will use Lemma \ref{lem:relationqes} and establish that $\beta_2(n)$ satisfies the formula announced for $(Y_n)$ in Theorem \ref{theo:dfreghrtu} . 
		
		Let $n$ points $U_1,\cdots,U_n$ with distribution $g_2$ in the triangle $\triangle$ with vertices $(0,0)$,$(1,0)$,$(0,3)$ as represented on Fig. \ref{fig:lafig}. 
		Consider $V=(X,H)$ the $U_i$ with the smallest height. The points $(U_i)$ together with $(1,0)$ are in convex position iff all the $U_i$ (except $V$) are in the colored zone represented in Fig. \ref{fig:lafig}, and, with $V$, are in convex position. The weight of the colored zone for $g_2$ is
		\[W(a,h)=\int_{h}^{\frac{h}{1-a}}   \int_{0}^{1-\frac{(1-a)}hy}  2x\mathrm{d}x\,\mathrm{d}y=\frac{a^3h}{3(1-a)}.\]
		We claim then that 
		\ben\label{eq:qghrt}  Y_n &=& \int_{0}^3\int_{0}^{1-h/3} 2an W(a,h)^{n-1} Y_{n-1}   \mathrm{d}a\,\mathrm{d}h\\
		&=& Y_{n-1} \frac{2}{n(3n-1)}= \frac{2^n}{n!}\prod_{j=1}^n \frac{1}{3j-1}.\een
		To finish the proof, we need to explain the following facts: the factor $n$ comes from the fact that any point can be at position $(a,h)$, the factor $W(a,h)^{n-1}$ accounts for the presence probability of all other points in the colored zone, and $\beta_{n-1}$ accounts for the probability that these points together with $V$ are in convex position. In order to prove that this last quantity is $\beta_{n-1}$, we need to remark the following: let us rescale the colored region by applying the linear map
		\[\Phi:(x,y)\to \l(\frac{x}{a}, 3\frac{1-a}{ah}(y-h)\r)\] 
		which bijectively maps the colored region onto the initial $\triangle$, and maps vertical lines onto vertical lines. 
		If $U$ is taken under the density $g_2$, and conditioned to be in the colored region, then $\Phi(U)$ has density $g_2$ on $\triangle$ (since this is the only density null on the axis $x=0$, and linear over $\triangle$, no computation is needed to make this deduction). 
	\end{proof} 
	\begin{rem}
		In order to obtain a similar lower bound as Theorem \ref{theo:dfreghrtu}  for $Q_{P}(n)$ the Prism $P=F\times[0,1]$ case (in the 3D case), we would need to compute $(\beta_1(n),n\geq 0)$, but we did not succeed. 
	\end{rem}

	\subsection{Bound in the tetrahedron case} 
	
	Consider a tetrahedron $T$ with vertices 
	\[A=(0,0,0),B=(1,0,0),C=(0,1,0),D=(0,0,6),\] and floor $F = \CH({A,B,C})$. The value of $Q_T(n)$ does not depend on the side lengths of the tetrahedron (as long as $T$ is taken non trivial), but note that here the floor area is $1/2$: the simplicity of the floor will help to get a tractable formula in the sequel. Note that for this reason, "the apex" of $T$, seen as a mountain, is $D$ and has a height of 6.

	The aim of this section is to prove Theorems \ref{theo:upp} and \ref{theo:dqgrh}, that give two (non trivial) sequences $(\ell_n)$ and $(u_n)$ such that
	\ben \ell_n \leq Q_T(n) \leq u_n, \textrm{~~for~~} n\geq 0.\een
	
	\subsubsection{Computation of an upper bound.}
	
	\begin{theo}\label{theo:upp} For any $n$, we have $Q_T(n)\leq u_n$, where $(u_n)$ satisfies the following recursion: $u_0=1$, and for $n\geq 1$,
		\ben \label{eq:fgto} 
		u_n = \frac{6}{(n+2)(n+1)n} \sum_{k=0}^{n-1} u_k u_{n-1-k},\een
		(so that, in particular $u_1=Q_T(1)=1$, $u_2=Q_T(2)=1/2$).
		The first terms from $u_0$ to $u_9$ are:
		\[\l[1,1,{\frac{1}{2}},{\frac{1}{5}},{\frac{7}{100}},{\frac{79}{3500}},{
			\frac{337}{49000}},{\frac{2069}{1029000}},{\frac{7033}{12348000}}\r]\]
		
	\end{theo}

	\begin{proof} 
		Let us call top of $T$,  the triangle $BCD$. 
		We define a floor to top projection, $\pi_t:ABC\to BCD$ parallel to the normal to the floor. Hence,  
		for $q=(x,y,0)$,
		\[\pi_t(q):=\pi_t(x,y,0)=(x,y,h(x,y)),\]
		where  $h(x,y)= 6\,(1-(x+y))$ is the height of the top, above $q$. 
		
		\begin{figure}[htbp]
			\centerline{\includegraphics{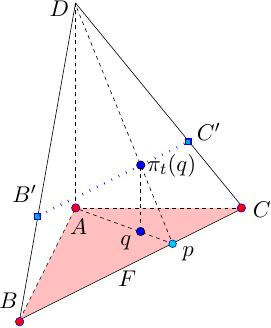}}
			\caption{\label{fig:pitq} Representation of $q$, and $\pi_t(q)$. The plane incident to $A,D,q,\pi_t(q)$ cuts the initial tetrahedron into two parts. }
		\end{figure}
		The vertical plane incident to $A,D,q,\pi_t(q)$ intersects the segment $(B,C)$ at the point \[p=\l(\frac{x}{x+y},\frac{y}{x+y},0\r),\]
		and this plane cuts the initial tetrahedron in two parts. In one of the two half-spaces, one finds the tetrahedron
		$T_1(q)=A,\pi_t(q),B,D$, and in the other side $T_2(q)=A,\pi_t(q),C,D$. 
		We have 
		\be \Vol_2(ABC)=\frac12, ~~
		\Vol_2(ABp)=\frac{y}{2(x+y)}, ~~
		\Vol_2(ACq)=\frac{x}{2(x+y)}.\ee Hence 
		\[\Vol_3(ADBp)=\frac{y}{x+y},~~~\Vol_3(ADCp)=\frac{x}{x+y},\] and then 
		\[\Vol_3(ABp\pi_t(q))=yh(x,y)/(6(x+y)), ~~\Vol_3(ACp\pi_t(q))=x h(x,y)/(6(x+y)).\] Finally set
		\be
		V_1(x,y):=\Vol_3(T_1(q))&=& \Vol_3(ADBp)-\Vol_3(ABp\pi_t(q))=\frac{y}{x+y}-y\frac{(1-(x+y))}{x+y}=y \\
		V_2(x,y):=\Vol_3(T_2(q))&=& \Vol_3(ADCp)-\Vol_3(ACp\pi_t(q))=\frac{x}{x+y}-x\frac{(1-(x+y))}{x+y}=x
		\ee
		The way of reasoning from here is very similar to that which led to "the comb formulas" in \cite{marckert:hal-02913348}.
		{We will need two main principles:\\
			(1) If $(X,Y,Z)$ is a uniform point in the tetrahedron, then conditional on $(X,Y)$, the variable $Z$ is uniform in $[0,h(X,Y)]$. To prove this, take a bounded continuous test function $f:\R^3\to \R$, and write $\E(f( Z/h(X,Y) , X,Y))=\int_{x+y\leq 1} \int_{z=0}^{h(x,y)} f(z/h(x,y),x,y) \d z \d x\d y$, then proceed to the change of variables $(u,x,y)=(z/h(x,y),x,y)$, so that $\d z\d x\d y=h(x,y)\d u\d x\d y$, and get $\E( f(Z/h(X,Y),X,Y))=\int_{x+y\leq 1} \int_{u=0}^{1} f(u,x,y) \d u (h(x,y) \d x\d y)$, from what we see that the marginal $Z/h(X,Y)$ is uniform and independent of $(X,Y)$, whose distribution has density $h(x,y)\1_{x+y\leq 1} \d x\d y$.\\
			(2) Take $(U_i,1\leq i \leq n)$  some iid uniform random variables taken in a set $D$ (with positive Lebesgue measure in $\R^d$), and $D'\subset D$, also with positive Lebesgue measure in $\R^d$. Then, if we condition by the event $E=\{ U_i \in D', 1\leq i\leq n\}$, then the $U_i$ are independent, and uniform in $D'$.}
		
		Now, consider the coordinates of the points $U_i=(X_i,Y_i,Z_i)$ for the independent uniform random points taken in $T$. Conditional on $(X_i,Y_i)$, $Z_i$ being uniform on $[0,h(X_i,Y_i)]$ can be written as $Z_i=V_i h(X_i,Y_i)$ with $V_i$ uniform on $[0,1]$, independent of $X_i,Y_i$. Now, consider $i^\star$ the special index $i$ of the maximal $V_i$. Conditional on $(i^\star,V_{i^\star})$, the other $V_i$ are uniform on $[0,V_{i^\star}]$ (and independent). 
		In other words, the other variables $U_i$ are iid uniform below the plane $(x,y,h(x,y)V_{i\star})$ while $U_{i^\star}$ is on that plane. It means that the other variables can be written under the form $(X_i,Y_i,h(X_i,Y_i)V_{i\star}W_i)$ where the $W_i$ are uniform on $[0,1]$, and the $(X_i,Y_i)$ are  distributed according to the density $h(x,y) \d x \d y$ as usual (and all these variables are independent).\\
		In words, the points $(X_j,Y_j,Z_j)$ are now written under the form $(X_j,Y_j,h(X_j,Y_j)W_jV_{i^\star})$ where the $W_j$ are uniform on $[0,1]$, except for $W_{i^\star}$ which is equal to 1.
		The fact that they are in convex position with the floor is equivalent to the problem where we remove $V_{i^\star}$.
		
		From here we can deduce that   
		\[Q_T(n)=Q'_T(1,n-1):=\P\l( \{U_0,\cdots,U_{n-1}\}\subset \partial CH(F \cup \{U_0,\cdots,U_{n-1}\})\r)\] where:\\
		-- the point $U_0$ is placed at position $(X,Y,h(X,Y))$ on the surface of the top triangle, where the density of the variable $(X,Y)$ is $h(x,y)\d x \d y$ (notice that $(X,Y)$ is not biased because of the independence between the $(V_i)$ and the $(X_i,Y_i)$),\\
		-- the points $U_1,\cdots,U_{n-1}$ are inside $T$ uniformly at random, independently from $U_0$. \\
		Conditionally on $U_0=(x,y,h(x,y))$ (which is $\pi_t(q)$ in the previous discussion), consider the subsets of the indices $J_1$ and $J_2$ of the  $U_i$ belonging to $T_1(q)$ and $T_2(q)$. If the points $U_0,U_1,\cdots,U_{n-1} \in \partial CH(F \cup \{U_0,\cdots,U_{n-1}\})$, then $J_1$ and $J_2$ must form a partition of $\{1,\cdots,n-1\}$. In particular, we must have
		\ben\label{eq:rge}
		\{U_i, i\in J_1\} &\in& \partial \CH(\{U_i, i\in J_1\} \cup \{A,B,\pi_t(q)\},\\
		\{U_i, i\in J_2\} &\in& \partial \CH(\{U_i, i\in J_2\} \cup \{A,C,\pi_t(q)\}
		\een so that, we are back in the initial floor problem settings.
		We then deduce that 
		\be
		Q_T(n)&=&Q_T'(1,n-1)\\
		&\leq& \int_{(x,y)\in ABC} \sum_{k=0}^{n-1}\binom{n-1}k  h(x,y) \Vol(T_1)^{k}\Vol(T_2)^{n-k-1} Q_T(k) Q_T(n-k-1)  dxdy.
		\ee
		The symbol "$\leq$" comes from the fact that the union over $(J_1,J_2)$ of the events \eref{eq:rge} contains positions of points that are not globally in convex position. 
		Observe that
		\be
		\alpha(n,k)
		&:=&\binom{n-1}k\int_{(x,y) \in ABC} V_1(x,y)^k V_2(x,y)^{n-1-k} h(x,y) dxdy\\
		&=&\binom{n-1}k\int_{x+y\leq 1\atop x,y\geq 0} y^k x^{n-1-k} 6(1-(x+y))dxdy
		\ee

		Using again \eref{eq:beta}, 
		\be
		D[k,j]:= \int_{x=0}^1 x^{k} \int_{y=0}^{1-x} y^j dydx     &=& \int_{x=0}^1 x^{k}  (1-x)^{j+1}/(j+1)dx  
		=  \frac{k! j!}{(k+j+2)!}
		\ee
		and then
		\[\alpha(n,k)= 6\binom{n-1}k (D[k,n-1-k]- D[k,n-k]-D[k+1,n-1-k])=6 \frac{(n-1)!}{(n+2)!}.\]
		In fine, putting all together,
		$Q_T(n) \leq \sum_{k=0}^{n-1}  \frac{6 (n-1)!} {(n+2)!}Q_T(k)Q_T(n-1-k)$ from what the result follows.
	\end{proof}

	\subsubsection{Lower bound}
	
	\begin{theo}\label{theo:dqgrh}  Let $(\ell_n,n\geq 0)$ be the sequence defined by $\ell_0=1$, for 
		\ben \ell_n= 6 \frac{(n-1)!\; n!}{(2n+1)!}\sum_{k=0}^{n-1} \ell_k\; \ell_{n-1-k},~~~\textrm{ for }~~n\geq 1,\een
		then $Q_T(n)\geq \ell_n$ for all $n$.
	\end{theo}
	The first terms of the sequence $(\ell_n)$, for $n$ from 0 to 6 are
	\be
	\left[1,1,{\frac{1}{5}},{\frac{1}{50}},{\frac{11}{10500}},{\frac{431}{
			12127500}},{\frac{2371}{2801452500}}\right].
	\ee
	This bound is better than the general bound given in Theorem \ref{theo:dfreghrtu} , meaning $Y_n\leq \ell_n$ for all $n\geq0$ (we prove this statement at the end of this proof).
	
	\begin{proof}
		
		We proceed in the beginning as in the proof  of Theorem \ref{theo:upp} until the addition of the point $U_0$, and again use that $Q_T(n)=Q'_T(1,n-1)$. 
		
		Consider $B'C'$ the segment obtained as the intersection of the top triangle $BCD$ with the line parallel to $BC$ containing $\pi_t(q)$, oriented as in Figure \ref{fig:pitq}.
		
		Now, take the indices $J_1$ and $J_2$ of the $U_i$ (with $1\leq i \leq n-1$), such that $i\in J_1 \iff U_i \in T'_1:=AB\pi_t(q)B'$ and $i\in J_2\iff U_i\in T'_2:=AC\pi_t(q)C'$.
		
		We claim that if $J_1$ and $J_2$ form a partition of $\{1,\cdots,n-1\}$ and if
		\be
		\{U_i,i\in J_1\} &\subset& \partial\CH( \{U_i,i\in J_1\}\cup AB\pi_t(q)) \\
		\{U_i,i\in J_2\} &\subset& \partial\CH( \{U_i,i\in J_2\}\cup AC\pi_t(q)) 
		\ee
		then, the points $U_0,\cdots,U_{n-1} \in \CH( \{U_0,\cdots,U_{n-1}\}\cup ABC) $. Indeed, the plane $AB'C'$ is a supporting plane of the $\CH(\{U_0,\cdots,U_{n-1}\} \cup\{A,B,C\})$ which, once more, splits the points $U_i$ in two independent regions.

		Conditional on $U_0=\pi_t(q)=(x,y,h(x,y))$, $B'=(x,0,h(x,y))$ and $C'=(0,y,h(x,y))$
		we get (using the determinant formula for the computation of the volume),
		\ben
		V'_1(x,y)&:=&\Vol(T_1')= y(1-(x+y)),\\
		V'_2(x,y)&:=&\Vol(T_2')=x(1-(x+y)).
		\een
		We then have that a lower bound is given by
		\[\ell_n=\sum_{k=0}^{n-1} \ell_k\;\ell_{n-1-k}\; g(n,k)\]
		with
		\be
		g(n,k)&:=&\binom{n-1}k \int_{ABC}V_1'^k(x,y) V_2'^{n-1-k}(x,y) h(x,y)dxdy\\
		&=&6\binom{n-1}k \int_{ABC} y^k x^{n-1-k} (1-(x+y))^{n}  dx dy
		\ee
		Set $t=x+y$, and since $x+y\leq 1$ we get by expanding $y^k=\sum_{k_1=0}^k\binom{k}{k_1}t^{k_1}(-x)^{k-k_1}$
		\be
		g(n,k)&:=&6\binom{n-1}k \int\int_{0\leq x \leq t\leq 1}  \sum_{k_1=0}^k\binom{k}{k_1}(-1)^{k-k_1}t^{k_1}x^{n-1-k_1} (1-t)^n dxdt\\
		&=&6\binom{n-1}k \int_{0 \leq t\leq 1}  \sum_{k_1=0}^k\binom{k}{k_1}(-1)^{k-k_1}\frac{1}{n-k_1}t^{n} (1-t)^n dt\\
		&=&6\binom{n-1}k \frac{(n!)^2}{(2n+1)!} \underbrace{\sum_{k_1=0}^k\binom{k}{k_1}(-1)^{k-k_1}\frac{1}{n-k_1}}_{S_n}
		\ee
		we claim that this last sum $S_n$ equals $(n \binom{n-1}{k})^{-1}$ so that we get $g(n,k):=6 n!(n-1)!/(2n+1)!$ in the end, which is the theorem statement. 
		
		To prove the claim, set $\psi(x):=\sum_{k_1=0}^k\binom{k}{k_1}(-1)^{k-k_1}\frac{x^{n-k_1}}{n-k_1}$ and observe that \[\psi'(x)= \sum_{k_1=0}^k\binom{k}{k_1}(-1)^{k-k_1} x^{(n-1)-k_1}=x^{n-1} (-1+1/x)^{k}=x^{n-1-k}(1-x)^{k}.\] Since $\psi$ is the primitive of $\psi'$ which is null at zero, we have $S_n=\psi(1)=\int_0^1\psi'(x)dx=\int_0^1x^{n-1-k}(1-x)^{k}dx $, so that with \eref{eq:beta}, $S_n=(n \binom{n-1}{k})^{-1}$.  
		
		To prove that $Y_n\leq \ell_n$ where $Y_n$ was introduced in  Theorem \ref{theo:dfreghrtu} , we proceed by induction and assume that $\ell_k\geq Y_k$ for all $1\leq k\leq n-1$ (one can check that this property is indeed true for the first values of $k$).
		Now, by definition of $\ell_n$ and $Y_n$, and by hypothesis, we have the inequality
		\begin{align*}
			\ell_n &\geq \frac{6}{n(2n+1)}{2n \choose n}^{-1}\sum_{k=0}^{n-1}Y_k Y_{n-1-k}\\
			&\geq \frac{6}{n(2n+1)}{2n \choose n}^{-1} \frac{2^{n-1}}{(n-1)!}\sum_{k=0}^{n-1}{n-1 \choose k} \l(\prod_{j=1}^{k}\frac{1}{3j-1}\r)\l(\prod_{j=1}^{n-1-k}\frac{1}{3j-1}\r)\\
			&\geq \frac{6}{n(2n+1)}{2n \choose n}^{-1} Y_{n-1}\sum_{k=0}^{n-1}{n-1 \choose k}\frac{\prod_{j=k+1}^{n-1}{(3j-1)}}{\prod_{j=1}^{n-1-k}(3j-1)}
		\end{align*}
		By factor comparisons, one can check that $\frac{\prod_{j=k+1}^{n-1}{(3j-1)}}{\prod_{j=1}^{n-1-k}(3j-1)}\geq {n-1 \choose k}$, so that, since $\binom{2(n-1)}{n-1}=\sum_{k=0}^{n-1}\binom{n-1}k^2$, we reach
		\[\ell_n \geq\frac{6}{n(2n+1)}\frac{{2(n-1)\choose n-1}}{{2n \choose n}}Y_{n-1},\]
		and it then suffices to see that $\displaystyle\frac{6}{n(2n+1)}\frac{{2(n-1)\choose n-1}}{{2n \choose n}}\geq \frac{2}{n(3n-1)}$ to complete the claim.
		
	\end{proof} 
	\subsection*{Acknowledgments}
	We would like to thank both reviewers for the precious time and efforts they spent reviewing this
	manuscript. We sincerely appreciate all their valuable comments and suggestions, which have helped
	us improve the quality of this paper.

	\normalsize

	\bibliographystyle{abbrv}

\begin{thebibliography}{10}
		
		\bibitem{barany1}
		I.~B{\'a}r{\'a}ny.
		\newblock Affine perimeter and limit shape.
		\newblock {\em Journal für die reine und angewandte Mathematik}, 484:71--84,
		1997.
		
		\bibitem{barany2}
		I.~B{\'a}r{\'a}ny.
		\newblock Sylvester's question: The probability that n points are in convex
		position.
		\newblock {\em The Annals of Probability}, 27(4):2020--2034, 1999.
		
		\bibitem{baranybipointed}
		I.~B{\'a}r{\'a}ny, G.~Rote, W.~Steiger, and C.-H. Zhang.
		\newblock A central limit theorem for convex chains in the square.
		\newblock {\em Discrete \& Computational Geometry}, 23:35--50, 01 2000.
		
		\bibitem{blaschke1917affine}
		W.~Blaschke.
		\newblock {\"U}ber affine geometrie xi: L{\"o}sung des “vierpunktproblems”
		von {S}ylvester aus der theorie der geometrischen wahrscheinlichkeiten.
		\newblock {\em Leipziger Berichte}, 69:436--453, 1917.
		
		\bibitem{Buchta1986}
		C.~Buchta.
		\newblock On a conjecture of {R}. {E}. {Miles} about the convex hull of random
		points.
		\newblock {\em Monatshefte für Mathematik}, 102:91--102, 1986.
		
		\bibitem{buchta_2006}
		C.~Buchta.
		\newblock The exact distribution of the number of vertices of a random convex
		chain.
		\newblock {\em Mathematika}, 53(2):247–254, 2006.
		
		\bibitem{buchta2012}
		C.~Buchta.
		\newblock On the number of vertices of the convex hull of random points in a
		square and a triangle.
		\newblock {\em Anzeiger. Abteilung II. Österreichische Akademie der
			Wissenschaften. Mathematisch-Naturwissenschaftliche Klasse}, 143, 01 2012.
		
		\bibitem{BR2001}
		C.~Buchta and M.~Reitzner.
		\newblock The convex hull of random points in a tetrahedron: solution of
		{B}laschke's problem and more general results.
		\newblock {\em J. Reine Angew. Math.}, 536:1--29, 2001.
		
		\bibitem{DM}
		A.~Dvoretzky and T.~Motzkin.
		\newblock A problem of arrangements.
		\newblock {\em Duke Math. J.}, 14:305--313, 1947.
		
		\bibitem{efronformula}
		B.~Efron.
		\newblock The convex hull of a random set of points.
		\newblock {\em Biometrika}, 52(3/4):331--343, 1965.
		
		\bibitem{groemer1974mean}
		H.~Groemer.
		\newblock On the mean value of the volume of a random polytope in a convex set.
		\newblock {\em Archiv der Mathematik}, 25(1):86--90, 1974.
		
		\bibitem{GT2021}
		A.~Gusakova and C.~Th\"ale.
		\newblock On random convex chains, orthogonal polynomials, {PF} sequences and
		probabilistic limit theorems.
		\newblock {\em Mathematika}, 67(2):434--446, 2021.
		
		\bibitem{Hostinsky}
		B.~Hostinsky.
		\newblock {\em Sur les probabilites geometriques}.
		\newblock Publ. Fac. Sci. Univ. Brno., 1925.
		
		\bibitem{kingman}
		J.~F.~C. Kingman.
		\newblock Random secants of a convex body.
		\newblock {\em Journal of Applied Probability}, 6(3):660--672, 1969.
		
		\bibitem{marckert2017probability}
		J.-F. Marckert.
		\newblock The probability that n random points in a disk are in convex
		position.
		\newblock {\em Brazilian Journal of Probability and Statistics},
		31(2):320--337, 2017.
		
		\bibitem{marckert:hal-02913348}
		J.-F. Marckert and S.~Rahmani.
		\newblock {Around Sylvester's question in the plane}.
		\newblock {\em {Mathematika}}, 67(4):860--884, Aug. 2021.
		
		\bibitem{miles}
		R.~E. Miles.
		\newblock Isotropic random simplices.
		\newblock {\em Advances in Applied Probability}, 3(2):353--382, 1971.
		
		\bibitem{morin2024bis}
		L.~Morin.
		\newblock Probability that $n$ points are in convex position in a general
		convex polygon: Asymptotic results.
		\newblock {\em {\it to appear in} Advances in Applied Probability}.
		
		\bibitem{morin2024}
		L.~Morin.
		\newblock Probability that {$n$} points are in convex position in a regular
		{$\kappa$}-gon: asymptotic results.
		\newblock {\em Adv. in Appl. Probab.}, 57(3):811--870, 2025.
		
		\bibitem{pfiefer}
		R.~E. Pfiefer.
		\newblock The historical development of {J. J. Sylvester's} four point problem.
		\newblock {\em Mathematics Magazine}, 62(5):309--317, 1989.
		
		\bibitem{MR2245368}
		J.~Pitman.
		\newblock {\em Combinatorial stochastic processes}, volume 1875 of {\em Lecture
			Notes in Mathematics}.
		\newblock Springer-Verlag, Berlin, 2006.
		\newblock Lectures from the 32nd Summer School on Probability Theory held in
		Saint-Flour, July 7--24, 2002, With a foreword by Jean Picard.
		
		\bibitem{MR1216521}
		R.~Schneider.
		\newblock {\em Convex bodies: the {B}runn-{M}inkowski theory}, volume~44 of
		{\em Encyclopedia of Mathematics and its Applications}.
		\newblock Cambridge University Press, Cambridge, 1993.
		
		\bibitem{sylvester}
		J.~J. Sylvester.
		\newblock Problem 1491.
		\newblock {\em The educational Times}, 1864.
		
		\bibitem{Valtr1995}
		P.~Valtr.
		\newblock Probability that n random points are in convex position.
		\newblock {\em Discrete and computational geometry}, 13(3-4):637--643, 1995.
		
		\bibitem{valtr1996probability}
		P.~Valtr.
		\newblock The probability that n random points in a triangle are in convex
		position.
		\newblock {\em Combinatorica}, 16(4):567--573, 1996.
		
	\end{thebibliography}

	\newpage
	\tableofcontents

\end{document}